\documentclass[12pt]{article}
\usepackage{amsmath}
\usepackage{latexsym}
\usepackage{amssymb}
\usepackage{a4wide}
\newtheorem{thm}{Theorem}[section]
\newtheorem{la}[thm]{Lemma}
\newtheorem{Defn}[thm]{Definition}
\newtheorem{Remark}[thm]{Remark}
\newtheorem{Note}[thm]{Note}
\newtheorem{prop}[thm]{Proposition}

\newtheorem{cor}[thm]{Corollary}
\newtheorem{Example}[thm]{Example}
\newtheorem{Examples}[thm]{Examples}
\newtheorem{Problems}[thm]{Problems}

\newtheorem{Problem}[thm]{Problem}
\newtheorem{Number}[thm]{\!\!}
\newenvironment{defn}{\begin{Defn}\rm}{\end{Defn}}

\newenvironment{example}{\begin{Example}\rm}{\end{Example}}

\newenvironment{rem}{\begin{Remark}\rm}{\end{Remark}}
\newenvironment{numba}{\begin{Number}\rm}{\end{Number}}
\newenvironment{proof}{{\noindent\bf Proof.}}%
                  {\nopagebreak\hspace*{\fill}$\Box$\medskip\smallskip\par}   
\newcommand{\Punkt}{\nopagebreak\hspace*{\fill}$\Box$}
\newcommand{\wb}{\overline}

\newcommand{\wt}{\widetilde}

\newcommand{\n}{\rm}

\newcommand{\impl}{\Rightarrow}
\newcommand{\ident}{\equiv}
\newcommand{\mto}{\mapsto}

\newcommand{\ve}{\varepsilon}
\newcommand{\isom}{\cong}

\newcommand{\Ad}{\mbox{{\rm Ad}}}

\newcommand{\N}{{\mathbb N}}

\newcommand{\R}{{\mathbb R}}

\newcommand{\Q}{{\mathbb Q}}

\newcommand{\Z}{{\mathbb Z}}

\newcommand{\Primes}{{\mathbb P}}
\newcommand{\Mix}{{\mathbb {MIX}}}
\newcommand{\Sub}{{\mathbb {SUB}}}
\newcommand{\LSub}{{\mathbb {VSUB}}}

\newcommand{\Lie}{{\mathbb {LIE}}}
\newcommand{\Loc}{{\mathbb {LOC}}}

\newcommand{\K}{{\mathbb K}}

\newcommand{\A}{{\mathbb A}}

\newcommand{\cV}{{\cal V}}

\newcommand{\g}{{\frak g}}
\newcommand{\ck}{{\frak k}}

\newcommand{\cp}{{\mathfrak p}}
\newcommand{\cq}{{\mathfrak q}}
\newcommand{\Cr}{{\mathfrak r}}
\newcommand{\dl}{{\displaystyle \lim_{\longrightarrow}}}

\newcommand{\Aut}{\mbox{\n Aut}}

\newcommand{\one}{\mbox{\rm \bf 1}}

\newcommand{\sub}{\subseteq}
\newcommand{\dt}{\mbox{\rm det}}

\newcommand{\SL}{\mbox{\rm SL}}

\newcommand{\im}{\mbox{\n im}}

\newcommand{\car}{\mbox{{\rm char}}}

\newcommand{\pr}{\mbox{\rm pr}}

\newcommand{\id}{\mbox{\n id}}

\newcommand{\bL}{{\mathbb L}}

\newcommand{\spr}{\mbox{\n\footnotesize pr}}

\begin{document}
\begin{center}
{\Large \bf Locally compact groups built up from
{\boldmath $p$\/}-adic\\[1.7mm]
Lie groups,
for {\boldmath $p$\/} in a given set of primes}\vspace{4.1 mm}\\
{\bf Helge Gl\"{o}ckner}\vspace{.7mm}
\end{center}
\noindent{\bf Abstract.\/}
We analyze the structure of locally
compact groups which can be built up
from $p$-adic Lie groups,
for $p$ in a given set of primes.
In particular, we calculate the scale function
and determine tidy subgroups for such groups,
and use them to recover the primes needed
to build up the group.\\[3mm]
{\footnotesize
}\vspace{3mm}
\begin{center}
{\large\bf Introduction}
\end{center}
While connected locally compact groups
can be approximated by real Lie groups
and hence can be described using real Lie theory
(\cite{MaZ}, \cite{HaM}),
the situation is more complicated
in the case of a totally disconnected,
locally compact group~$G$.
Here,
we have a $p$-adic Lie theory
available for each prime~$p$,
but it is not clear {\em a priori\/}
which primes~$p$
will be needed to analyze the structure
of~$G$, nor whether $p$-adic Lie theory
is useful at all in this context.
Investigations in \cite{App}
indicate that indeed general locally compact groups
are ``too far away'' from $p$-adic Lie groups
(and from Lie groups over local fields) to expect
meaningful applications of Lie theory.
Therefore, it is essential to
restrict attention
to suitable classes of totally disconnected groups,
which are ``close enough'' to $p$-adic Lie groups.\\[3mm]
For example, we might consider the class
of (locally compact) pro-$p$-adic Lie groups,
viz.\ locally compact groups which can be
approximated by $p$-adic Lie groups, for a fixed prime~$p$
(see \cite{App}, \cite{GaW}
for investigations of
such groups). However, it is clearly
very restrictive to use $p$-adic
Lie theory for a single prime~$p$ only;
it would be more
natural to try to
make use of $p$-adic Lie theory for variable
primes $p$ simultaneously.
For instance, it should certainly be allowed to
approximate a group also by finite products
$\prod_{p\in \cp}G_p$
of $p$-adic Lie groups~$G_p$ (where $\cp$ is a finite set of primes),
or by closed subgroups of
such products.\\[3mm]
Motivated by such considerations, given a non-empty subset
$\cp$ of the set $\Primes$ of all primes,
it was proposed
in \cite{App} to study the class $\Mix_\cp$ of all
locally compact groups which can be manufactured
from $p$-adic Lie groups with $p\in \cp$,
by repeated application of the operations
of forming cartesian products,
closed subgroups, Hausdorff quotients,
and passage to isomorphic topological groups.\footnote{Thus,
technically speaking, $\Mix_\cp$ consists of all
locally compact groups in the variety of Hausdorff
groups generated by the class of topological groups
which are $p$-adic Lie groups for some $p\in \cp$.}\\[3mm]
For example, consider a Hausdorff
quotient $S/N$ of a closed subgroup~$S$
of a finite product
$\prod_{p\in F}G_p$
of $p$-adic Lie groups,
with $p$ in a finite set $F\sub \cp$,
or a topological group isomorphic to $S/N$ (such groups
will be called {\em $\A_\cp$-groups\/}).
Then $S/N$ is a $\Mix_\cp$-group.
Arbitrary $\Mix_\cp$-groups are not too far away from
this example: a locally compact group $G$ is a $\Mix_\cp$-group
if and only if it is topologically isomorphic to a closed
subgroup of a cartesian product
$\prod_{i\in I}S_i/N_i$ of $\A_\cp$-groups (cf.\ (\ref{subprod})),
by standard facts from the theory of varieties of
topological groups (\cite{BMS}, \cite{HMS},
\cite{Mor}).\\[3mm]
This information alone would
not be enough to analyze $\Mix_\cp$-groups
via $p$-adic Lie theory. However,
we can prove much more:
Every $G\in \Mix_\cp$
can be approximated by $\A_\cp$-groups,
in the sense that every identity neighbourhood
of $G$ contains a closed normal subgroup
$K\sub G$ such that $G/K$ is an $\A_\cp$-group (Remark~\ref{charlsub}).
We can also show that every $\A_\cp$-group~$G$
contains an open subgroup
which is a finite product $\prod_{p\in F}H_p$
of $p$-adic Lie groups (Corollary~\ref{Apthensub}).
Since every
inner automorphism of~$G$
gives rise to local automorphisms
of the factors $H_p$ here, adapting techniques from
\cite{Sca} and \cite{GaW}
to the case of local automorphisms
we are able to deduce
very satisfactory results concerning
the structure of $\Mix_\cp$-groups
(including solutions to all open
problems formulated in~\cite{App}).
In particular, we obtain
a clear picture of
the ``tidy subgroups'' of a $\Mix_\cp$-group~$G$
and its ``scale function'' $s_G\!: G\to \N$,
which are the essential structural features
of~$G$ in the
structure theory of
totally disconnected groups
initiated in (\cite{Wi1}, \cite{Wi3}).
We recall the definitions:\\[3mm]
{\bf Definition}
(cf.\ \cite{Wi1}, \cite{Wi3}).
Let $G$ be a
totally disconnected, locally compact group
and $\alpha$ be an automorphism
of~$G$.
A compact, open subgroup $U$ of $G$ is called
{\em tidy for~$\alpha$\/}
if the following conditions are satisfied:
\begin{itemize}
\item[(T1)]
$U=U_+U_-$, where $U_\pm:=\bigcap_{n\in \N_0}\alpha^{\pm n}(U)$;
\item[(T2)]
The subgroup
$U_{++}:=\bigcup_{n\in \N_0}\alpha^n(U_+)$ is closed in~$G$.
\end{itemize}
It can be shown that
compact, open subgroups tidy for $\alpha$
always exist, and that the index
\[
r_G(\alpha)\; :=\; [\alpha(U_+):U_+]
\]
(called the ``scale of~$\alpha$'')
is finite and independent of the choice of tidy subgroup~$U$.
Specialization to inner automorphisms
$I_x\!: G\to G$, $I_x(y):=xyx^{-1}$ yields the
{\em scale function\/}
$s_G\!: G\to \N$, $s_G(x):=r_G(I_x)$ of~$G$.
We let $\Primes(G)$ be the set of all
primes $p\in \Primes$ such that $p$ divides $s_G(x)$
for some $x\in G$.\\[3mm]
{\bf The main results.}
Writing $\Mix_\emptyset$ for the class of locally compact
pro-discrete groups, we can summarize
our main results as follows:
\begin{itemize}
\item[(a)]
For any sets of primes $\cp$ and $\cq$, we have
$\Mix_\cp\cap \Mix_{\cq}=\Mix_{\cp\cap\cq}$
(Theorem~\ref{intersections}).
\item[(b)]
The scale function $s_G$ of any
$\Mix_\cp$-group~$G$ can be calculated
by Lie-theoretic methods.
Furthermore, a basis of compact, open subgroups
tidy for $x$ can be described explicitly
for each $x\in G$, using Lie-theoretic
methods (Theorem~\ref{scaletidy},
Corollary~\ref{scalelim}).\footnote{This is new
even for $p$-adic groups;
in~\cite{Sca}, $s_G$ was calculated without
formulas for tidy subgroups.}
\item[(c)]
For every $\Mix_\cp$-group~$G$,
the set $\Primes(G)$ of all prime
divisors of the values of the scale function
is a finite set, and $\Primes(G)\sub \cp$
(Corollary~\ref{scalelim}).
\item[(d)]
If $G$ is a compactly generated $\Mix_\Primes$-group
and $\cp$ a set of primes, then $G\in \Mix_\cp$
if and only if $\Primes(G)\sub \cp$
(Theorem~\ref{exactprimes}).
In particular, finitely many primes suffice
to build up~$G$.
Furthermore, every compactly generated, uniscalar
$\Mix_\Primes$-group~$G$ is pro-discrete
(Corollary~\ref{unispro}).\footnote{Recall that
a totally disconnected,
locally compact group~$G$ is called {\em uniscalar\/} if $s_G\ident 1$,
which holds if and only if every $x\in G$
normalizes some compact, open subgroup of~$G$.}
Previously, this was only known
for $p$-adic Lie groups
(see \cite{Par} and \cite{GaW}).
\end{itemize}
It is a natural idea
that the set $\Primes(G)$
of all prime divisors of the values of $s_G$
should tell us which kinds
of $p$-adic Lie groups (which $p$) are needed to analyze
a totally disconnected group~$G$,
at least in good cases.\footnote{This idea was
expressed by M. Stroppel (Stuttgart)
in 1994.}
Result\,(d) above shows that
this general philosophy can be turned into a mathematical
fact for the class of compactly generated $\Mix_\Primes$-groups.\\[3mm]
We mention that most of the results carry over
to the properly larger class $\LSub_\cp$
(subsuming $\Mix_\cp$)
of all locally compact groups
in the variety of Hausdorff groups
generated by
topological
groups having a direct product
$\prod_{p\in F}H_p$
of $p$-adic Lie groups as an open subgroup,
for $p$ in a finite subset
$F\sub \cp$. We therefore discuss such groups
in parallel.\\[3mm]
Although our studies
may remind the reader of Ad\`{e}le groups,
closer inspection shows
that the latter need not belong
to $\Mix_\Primes$, nor $\LSub_\Primes$
(see Remark~\ref{adele}).\\[3mm]
{\bf Variants.}
Some results
remain valid if $p$-adic Lie groups
are replaced by {\em locally pro-$p$ groups\/}
(groups possessing
a pro-$p$-group
as a compact open subgroup):
see Section~\ref{seclocpro}.\\[3mm]
{\bf Further results.}
Motivated by results in~\cite{Wi2},
in the final Section~\ref{seclocal}
we associate a set of primes
$\bL(G)$ to each totally disconnected,
locally compact group~$G$,
which only depends on the local isomorphism
type of~$G$
(the ``local prime content of~$G$'').
Since $\bL(G)$ contains all prime
divisors of the scale function,
it provides a means to
deduce information concerning the global structure
of~$G$ (its scale function)
from the local structure of~$G$.
Using the local prime content,
we show that for each $G\in \Mix_\Primes$,
there exists a unique
smallest set of primes~$\cp$
such that $G\in \Mix_\cp$
(Theorem~\ref{noncomp}, Remark~\ref{thesmallest}).
If $G$ is compactly
generated, then simply $\cp=\Primes(G)$,
as mentioned before.
If $G$ is not compactly
generated, then $\cp\not=\Primes(G)$
in general.
In this case, $\cp$
can still be determined in principle
(it is the ``intermediate prime
content'' of $G$, defined below),
but it is a less tangible invariant.\\[3mm]
In an appendix, which is of independent
interest, we describe topological groups~$G$ whose
normal subgroups~$N$ with $G/N$ a real (resp. $p$-adic)
Lie group do not form a filter basis.\\[3mm]
The present paper
uses (and generalizes) results and techniques
from \cite{Sca}, \cite{App}, \cite{GaW}
and~\cite{Wi2}.\vspace{1.5mm}
\section{Preliminaries and notation}
\begin{numba}\label{defnvariety}
Given a class of topological Hausdorff
groups $\Omega$, the variety of Hausdorff groups
generated by~$\Omega$
is the smallest class $\cV(\Omega)$
of Hausdorff groups containing $\Omega$
and closed under the operations of formation
of cartesian products ``${\rm C}$,''
subgroups ``${\rm S}$,''
Hausdorff quotients ``${\rm \wb{Q}}$,''
and passage to isomorphic topological groups
(which is understood and suppressed in the notation).
It is easy to see that
\begin{equation}\label{var1}
\cV(\Omega)\; =\; {\rm \wb{Q}SC}(\Omega)
\end{equation}
here (cf.\ \cite[Thm.\,1]{BMS} or \cite[Thm.\,6]{Mor}),
and it can be shown with more effort that
\begin{equation}\label{var2}
\cV(\Omega)\; =\; {\rm SC\wb{Q}\,\wb{S}P}(\Omega)
\end{equation}
(see \cite[Thm.\,2]{BMS},
or \cite[Thm.\,7]{Mor}),
where ``${\rm P}$'' denotes the formation of all
finite cartesian products,
and ``${\rm \wb{S}}$'' denotes formation of closed
subgroups (or isomorphic copies thereof, as above).
It is easy to see (cf.\ (\ref{var1}) above) that the class
${\rm \wb{Q}\,\wb{S}P}(\Omega)$ is
closed under the formation of finite cartesian products,
closed subgroups and Hausdorff quotients.
\end{numba}
\begin{numba}\label{defnmixtures}
Throughout the following,
$\Primes$ denotes the set of all primes.
Given $p\in \Primes$, we let $\Lie_p$
be the class of $p$-adic Lie groups;
given a non-empty subset $\cp\sub \Primes$, we set
$\Lie_\cp:=\bigcup_{p\in \cp}\Lie_p$.
According
to (\ref{var2}), the variety of Hausdorff groups
generated by $\Lie_\cp$ is given by
\begin{equation}\label{subprod}
\cV(\Lie_\cp)\;=\; {\rm SC}(\A_\cp)\,,
\qquad \mbox{where}\qquad
\A_\cp\;:=\; {\rm \wb{Q}\,\wb{S}P}(\Lie_\cp)
\end{equation}
is the class of all topological groups
isomorphic to a Hausdorff quotient $S/N$
of a closed subgroup $S$ of
a product $\prod_{p\in F}G_p$,
where $F\sub\cp$ is a finite subset
and $G_p$ a $p$-adic Lie group for each $p\in F$.
For later use, we let $\A_\emptyset$
be the class of discrete groups. We define
\[
\Mix_\cp\; :=\; \{G\in \cV(\Lie_\cp)\!: \,\mbox{$G$ is locally compact}\,\}\,.
\]
Finally, we let
$\Mix_\emptyset$ be the class of all {\em pro-discrete\/},
locally compact groups~$G$,
{\em i.e.},
locally compact groups~$G$
whose filter of identity neighbourhoods
has a basis consisting of open, normal subgroups of~$G$
(see \cite{App} for more information).
\end{numba}
\begin{numba}
Given a set $\cp$ of primes,
we let $\Sub_\cp$
be the class of all topological groups
possessing an open subgroup
isomorphic to $\prod_{p\in F}G_p$,
where $F\sub \cp$ is finite
and $G_p$ a $p$-adic Lie group
for each $p\in F$.
We let
\[
\LSub_\cp \; :=\;
\{G\in \cV(\Sub_\cp)\!: \,\mbox{$G$ is locally compact}\,\}\,.
\]
In particular,
$\LSub_\emptyset$ is the class of locally compact, pro-discrete
groups (cf.\ \cite[Thm.\,2.1]{App}).
\end{numba}
Note that all of the topological groups
in $\A_\cp$, $\Sub_\cp$,
$\Mix_\cp$ and $\LSub_\cp$
are locally compact and totally disconnected.
We shall see later that $\A_\cp\sub \Sub_\cp$
and thus $\Mix_\cp\sub \LSub_\cp$.
\begin{numba}\label{defnadm}
A class $\A$ of Hausdorff topological groups
which contains the trivial group and is closed under
passage to isomorphic topological groups is called
a {\em property of topological groups};
the elements of $\A$ are called {\em $\A$-groups}.
If $\A$ is a property of topological groups,
we say that the class $\A$
is {\em suitable for approximation\/}
(or also: an ``admissible property'' of topological groups,
in the terminology of \cite{App}),
if every $\A$-group is locally compact,
$\A$ is closed under the formation of
finite cartesian products,
closed subgroups and Hausdorff quotients
(which
holds if and only if $\A={\rm {\wb Q}\,{\wb S}P}(\A)$),
and if $G/{\ker f}$ is an $\A$-group,
for every continuous homomorphism $f\!: G\to H$
from a locally compact group~$G$ to an $\A$-group~$H$.
\end{numba}
For example, the class of real
Lie groups is suitable for approximation
(cf.\ \cite{HMS}),
and so are the classes
of $p$-adic Lie groups (see \cite{App}),
finite groups, finite $p$-groups,
and finite nilpotent groups,
respectively.\\[3mm]
Quite a bit of work will be needed to see
that the classes $\A_\cp$ and $\Sub_\cp$ are
suitable for approximation. This information
is very useful, because it is well understood
which locally compact groups can be approximated
by topological groups in a class of topological groups
which is suitable for approximation.
We recall \cite[Thm.\,2.1]{App}:
\begin{prop}\label{appradm}
Let $\A$ be a class of topological groups
that is suitable for approximation,
and $G$ be a locally compact group. Then the following
conditions are equivalent:
\begin{itemize}
\item[\rm (a)]
$G$ can be approximated by $\A$-groups, i.e.,
every identity neighbourhood~$U$ of~$\,G$
contains a closed normal subgroup $N$ of $\,G$
such that $G/N\in \A$.
\item[\rm (b)]
The set of all closed
normal subgroups $N$
as in {\rm (a)} is a filter basis which
converges to~$1$ in~$G$.
\item[\rm (c)]
$G$ is a pro-$\A$-group in the sense of {\rm \,\cite[1.4]{App}}.
\item[\rm (d)]
$G$ is a projective limit $($in the category of topological groups$)$
of a projective system of $\A$-groups and continuous
homomorphisms.
\item[\rm (e)]
$G$ is an element of the variety $\cV(\A)$ of Hausdorff
groups generated by~$\A$.\Punkt
\end{itemize}
\end{prop}
Choosing $U$ compact,
we see that $N$ in~(a) can always be assumed to be compact.
\begin{numba}
All topological groups considered in this article
are Hausdorff.
Open, surjective, continuous homomorphisms
are called {\em quotient morphisms}.
All isomorphisms or automorphisms
of topological groups are, in particular, homeomorphisms.
The automorphism group
of a topological group~$G$ is denoted $\Aut(G)$.
A {\em local isomorphism\/} between
totally disconnected,
locally compact groups
$G$ and $H$
is an isomorphism from an open subgroup of~$G$
onto an open subgroup of~$H$.
\end{numba}
\begin{numba}
Our main sources for $p$-adic Lie
theory are \cite{BLi} and \cite{Ser}.
All Lie groups~$G$ considered here
are finite-dimensional analytic Lie groups
(unless we say otherwise explicitly).
As usual, a $p$-adic Lie group
will be identified with its
underlying topological group.
The $p$-adic Lie algebra of~$G$
is denoted~$L(G)$.
All necessary background concerning
pro-finite groups
and pro-$p$-groups
(in particular, the basics
of Sylow theory needed here)
can be found in~\cite{Wls}.
\end{numba}
\begin{numba}\label{deflocpro}
Given a prime~$p$,
a topological group~$G$ is called
{\em locally pro-$p$\/}
if $G$ has a compact, open subgroup $U$
which is a pro-$p$-group.
As an immediate consequence of the corresponding
permanence properties of pro-$p$-groups,
the class $\Loc_p$ of locally pro-$p$ groups
is closed under formation of finite direct products,
closed subgroups, and Hausdorff quotients.
\end{numba}
\begin{numba}\label{wellkn}
It is well known
that
every $p$-adic Lie group
is locally pro-$p$,
and so is
every analytic Lie group~$G$ over a local
field~$\K$ whose
residue field~${\frak k}$ has
characteristic~$p$~\cite{Ser}.
See also \cite[Prop.\,2.1\,(h)]{ANA}
for a recent proof, which remains
valid if~$G$
is not analytic
but merely a $C^1$-Lie group
(in the setting of~\cite{BGN}).
While every $p$-adic $C^k$-Lie group
admits a $C^k$-compatible
analytic Lie group structure~\cite{ANA},
for every local field
of positive characteristic
there exists a $1$-dimensional
smooth Lie group without an
analytic Lie group structure
compatible with its topological
group structure,
and $C^k$-Lie groups which are not
$C^{k+1}$~\cite{NOA}.
\end{numba}
\section{Relations between the various classes of groups}
We first collect various simple facts.
\begin{la}\label{propproq}
Let $p\not=q$ be primes,
$G$ be a pro-$p$-group,
$H$ a pro-$q$-group
and $f\!: G\to H$
be a continuous
homomorphism.
Then $f(x)=1$
for all $x\in G$.
In particular, every continuous
homomorphism from a $p$-adic Lie group
to a $q$-adic Lie group
has open kernel.
\end{la}
\begin{proof}
Since continuous homomorphisms
to finite $q$-groups separate points on~$H$,
we may assume that~$H$ is a finite
$q$-group. By \cite[La.\,1.2.6]{Wls},
$K:=\ker \, f$ is open
and hence $G/K\isom f(G)$
is a $p$-group
(as a consequence of \cite[Prop.\,1.2.1]{Wls}).
Hence $f(G)=\{1\}$. The well-known
final assertion (Cartan's Theorem)
now follows with {\bf \ref{wellkn}}.
\end{proof}
The following observation
concerning closed subgroups of
pro-nilpotent groups is the key
to an understanding
of $\Mix_\cp$-groups and
$\LSub_\cp$-groups.\footnote{Recall from
\cite{Wls} that
projective limits
of nilpotent finite groups
are called {\em pro-nilpotent}.}
\begin{prop}\label{block}
Let
$U_p$ be a pro-$p$-group
for each $p\in \Primes$,
and $S$ be a closed subgroup
of $U:=\prod_{p\in \Primes}U_p$.
Define $S_p:=S\cap U_p$
for $p\in \Primes$, identifying $U_p$
with $U_p\times \prod_{q\in \Primes\setminus\{p\}}\{1\}\sub U$.
Then
$S=\prod_{p\in \Primes} S_p$.
\end{prop}
\begin{proof}
Being a direct product
of pro-$p$-groups, $U$ is
pro-nilpotent (cf.\
\cite[Prop.\,2.4.3]{Wls}).
As a consequence of
\cite[Thm.\,1.2.3]{Wls},
also the closed subgroup~$S$
of~$U$ is pro-nilpotent.
Hence~$S$ has a normal (and hence unique)
$p$-Sylow subgroup $\tilde{S}_p$
for each $p\in \Primes$
(see \cite[Prop.\,2.4.3\,(ii) and Prop.\,2.2.2\,(d)]{Wls}).
Since $\tilde{S}_p$
is contained in the unique $p$-Sylow
subgroup~$U_p$ of~$U$
and contains $S_p$
(see \cite[Prop.\,2.2.2\,(c)]{Wls}),
we deduce that
$\tilde{S}_p=S_p$.
Hence $S=\prod_{p\in \Primes} S_p$
by \cite[Prop.\,2.4.3\,(iii)]{Wls}.
\end{proof}
\begin{cor}\label{nicesub}
Let $\cp$ be a non-empty, finite set of primes,
$G_p$ be a locally pro-$p$ group
$($resp., a $p$-adic Lie group$)$ for $p\in \cp$,
and $S$ be a closed subgroup of $\,G:=\prod_{p\in \cp}G_p$,
which we consider as an internal direct product
of the groups $G_p$.
Then
$S_p:=S\cap G_p$ is a closed normal
subgroup of $S$
and locally pro-$p$ $($resp.,
a $p$-adic Lie group$)$,
being a closed subgroup of~$G_p$.
Furthermore, the product\vspace{-3mm}
\[
P\; := \; \prod_{p\in \cp}S_p\vspace{-5mm}
\]
is an open subgroup of~$S$.
\end{cor}
\begin{proof}
For each $p\in \cp$,
there exists an open pro-$p$ subgroup
$U_p\sub G_p$. Then $U:=\prod_{p\in \cp}U_p$
is open in~$G$ and hence
$U\cap S$ is open in~$S$,
where $U\cap S=\prod_{p\in \cp} (U_p\cap S)$
by Proposition~\ref{block}.
Since $U_p\cap S\sub S_p$,
also $P$ is open in~$S$.
\end{proof}
We now focus on $p$-adic groups.
Analogues for locally pro-$p$ groups
are outlined in Section~\ref{seclocpro}.
\begin{cor}\label{injthensub}
Suppose that
$\phi\!: G\to \prod_{p\in \cp}G_p$
is a continuous, injective homomorphism
from a locally compact group~$G$ into a product
of $p$-adic Lie groups $G_p$, for $p$ in some finite set
of primes~$\cp$. Then $G$ is a $\Sub_\cp$-group.
\end{cor}
\begin{proof}
Since $\phi$ is injective, $G$ is totally disconnected.
We choose a compact, open subgroup
$U\sub G$; then $U$ is isomorphic
to the closed subgroup $S:=\phi(U)$
of $\prod_{p\in \cp}G_p$.
Thus
Corollary~\ref{nicesub} entails the claim.
\end{proof}
\begin{cor}\label{Apthensub}
Let $\cp$ be a
set of primes.
Then $\A_\cp\sub \Sub_\cp$ and thus $\Mix_\cp\sub
\LSub_\cp$.
\end{cor}
\begin{proof}
Without loss of generality $\cp\not=\emptyset$,
the omitted case being trivial.
If $G$ is an $\A_\cp$-group,
then
after passing to an isomorphic copy we may
assume that $G=S/N$
where $S$ is a closed subgroup of a product
$P:=\prod_{p\in F}G_p$ of $p$-adic Lie groups
for $p$ in some finite subset $F\sub \cp$, and $N\sub S$
a closed normal subgroup.
Then $S_p:=S\cap G_p$ is a closed subgroup
of~$S$, and $N_p:=N\cap G_p$ is a closed normal subgroup
of~$S_p$, for each $p\in F$.
By Corollary~\ref{nicesub},
$\wt{S}:=\prod_{p\in F} S_p$ is an open subgroup of~$S$
and $\wt{N}:=\prod_{p\in F}\ N_p$ an open subgroup
of~$N$ (and it also is a closed normal subgroup
of~$S$). Hence
$G=S/N$ is isomorphic to
\begin{equation}\label{double}
\bigl(S / \wt{N}\bigr) \Big/ \bigl(N / \wt{N}\bigr)\;.
\end{equation}
But,
$N/\wt{N}$ being discrete, the
group in~(\ref{double})
is locally isomorphic to $S/\wt{N}$,
which has $\wt{S}/\wt{N}\isom \prod_{p\in F} \bigl(S_p/N_p\bigr)$
as an open subgroup. Hence for all
sufficiently small compact, open subgroups
$U_p\sub S_p/N_p$,
the product $\prod_{p\in F} U_p$
is isomorphic to a compact, open
subgroup of~$G$. Thus $G$ is a $\Sub_\cp$-group.
Hence $\A_\cp\sub \Sub_\cp$.
The rest is obvious.\vspace{-.6mm}
\end{proof}
%
%
%
%
%
We may assume that $\pr_p(S)$ is dense
in~$G_p$
in the preceding proof
(where $\pr_p\!: P\to G_p$
is the coordinate projection),
entailing that the closed
subgroup $N_p:=N\cap G_p$ of~$G_p$
is a normal subgroup of~$G_p$.
Hence $S/\wt{N}$ is a closed subgroup of $(\prod_{p\in F}G_p)/\wt{N}
\isom \prod_{p\in F}(G_p/N_p)$, where $G_p/N_p$ is a $p$-adic
Lie group. Combining this with (\ref{double}), we get:
\begin{cor}\label{additinfo}
If $\cp\not=\emptyset$,
then
every $\A_\cp$-group
is topologically isomorphic to a quotient
$S/D$, where $S$ is a closed
subgroup of a product
$\prod_{p\in F}G_p$
of $p$-adic Lie groups~$G_p$
for $p$ in a finite set $F\sub \cp$,
and $D$ is a {\bf discrete\/} normal subgroup
of~$S$.\Punkt
\end{cor}
\begin{prop}\label{admissible}
Let $\cp$ be a set of primes.
Then we have:
\begin{itemize}
\item[\rm (a)]
The class $\Sub_\cp$ is suitable for approximation.
\item[\rm (b)]
The class $\A_\cp$ is suitable for approximation.
\end{itemize}
\end{prop}
\begin{proof}
The case $\cp=\emptyset$ being trivial, we may assume that
$\cp$ is non-empty.

(a) Every $\Sub_\cp$-group is locally compact.
Let us show that the class $\Sub_\cp$
is closed under the formation of finite direct products,
closed subgroups, and Hausdorff quotients.
It is obvious that
finite products of $\Sub_\cp$-groups
are $\Sub_\cp$-groups.
Let $G$ be a $\Sub_\cp$-group now
and $S$ a closed subgroup.
Let $U\sub G$ be an open subgroup
which is a product $\prod_{p\in F}U_p$
of $p$-adic Lie groups $U_p$ for
$p$ in a finite subset $F\sub \cp$.
Then $S\cap U$ is a closed subgroup
of $U=\prod_{p\in F} U_p$,
whence $S$ is a
$\Sub_F$-group
(and thus {\em a fortiori\/}
a $\Sub_\cp$-group),
by Corollary~\ref{nicesub}.
If $G$ and $U$ are as before and $N$ is a closed normal subgroup
of~$G$, then $G/N$ has an open subgroup
isomorphic to
$U/(U\cap N)$, which is an $\A_\cp$-group and hence
a $\Sub_\cp$-group
by Corollary~\ref{Apthensub}.
Hence also $G/N$ is a $\Sub_\cp$-group.

Finally, suppose that $f\!: G\to H$ is a continuous
homomorphism from a locally compact group~$G$ to a
$\Sub_\cp$-group $H$. Since $H$ is totally disconnected,
the connected identity component $G_0$ of~$G$
is contained in the kernel of~$f$,
entailing that $Q:=G/\ker f$
is totally disconnected.
Let $\wb{f}\!: Q\to H$ be the injective continuous
homomorphism determined by $\wb{f}\circ q=f$,
where $q\!: G\to Q$ is the quotient map.
Let $W\sub H$ be an open subgroup
which is a finite product of $p$-adic Lie groups
(for certain $p\in \cp$).
Being totally disconnected and locally compact, $Q$
has a compact, open subgroup~$U$ contained in $\wb{f}^{\,{-1}}(W)$.
Then $U$ is a $\Sub_\cp$-group
by Corollary~\ref{injthensub},
and hence so is
$Q=G/\ker f$.
The proof of\,(a) is complete.

(b) It is obvious that $\A_\cp$ is closed under
the formation of closed subgroups,
Hausdorff quotients, and finite cartesian products.
Given a continuous homomorphism
$\phi\!: G\to H$ from a locally compact group
$G$ to an $\A_\cp$-group~$H$,
the locally compact group~$Q:=G/\ker\phi$
is totally disconnected. There is a unique continuous
injective homomorphism $\wb{\phi}\!: Q\to H$
such that $\wb{\phi}\circ \kappa=\phi$,
where $\kappa\!: G\to Q$ is the quotient map.
Let us show that~$Q$ is
an $\A_\cp$-group.
For convenience
of notation, after replacing $G$ with~$Q$
and $\phi$ with $\wb{\phi}$,
we may assume without loss
of generality that $\phi\!: G\to H$ is injective.
Furthermore,
in view of Corollary~\ref{additinfo},
we may assume that $H=S/D$
for some closed subgroup~$S$
of a product $\prod_{p\in F}G_p$
of $p$-adic Lie groups~$G_p$ for $p$
in a finite subset $F\sub \cp$,
and some discrete normal
subgroup~$D$ of~$S$. We let $\rho\!: S\to S/D=H$
be the canonical quotient map.
Our goal is to equip
$S':=\rho^{-1}(\phi(G))$
with a finer topology which turns this group
into a closed subgroup
of another, suitable chosen product of $p$-adic Lie groups,
and such that $S'/D\isom G$ as a topological group.

To this end,
we choose a compact, open subgroup
$U$ of~$G$.
Then
$W:=\rho^{-1}(\phi(U))$
is a closed subgroup of~$S$
and hence also of $\prod_{p\in F}G_p$.
Let $W_p:=W\cap G_p$;
by Corollary~\ref{nicesub},
$\prod_{p\in F}W_p$
is an open, normal subgroup of~$W$.
We claim that, for each $y\in S'$,
the inner automorphism $I_y\!: S\to S$,
$I_y(s):=ysy^{-1}$
restricts to
a local
automorphism of~$W$. Indeed,
there exists
$x\in G$ such that $\rho(y)=\phi(x)$,
and a compact, open subgroup $A\sub G$
such that $xAx^{-1}\sub U$.
Then, $\phi|_U^{\phi(U)}$
being a homeomorphism as~$U$ is compact,
$\phi(A)$ is an open, compact subgroup
of $\phi(U)$, such that $I_{\phi(x)}(\phi(A))=\phi(I_x(A))\sub \phi(U)$.
Hence, $\rho|_W^{\phi(U)}$
being continuous, $B:=\rho^{-1}(\phi(A))$
is an open subgroup of~$W$,
and from $\rho(I_y(B))=I_{\rho(y)}(\rho(B))=I_{\phi(x)}(\phi(A))
\sub \phi(U)$
we deduce that $I_y(B)\sub W$.
The continuity of $I_y\!: S\to S$
and its inverse
entails that $I_y|_B^W$ is an isomorphism
from~$B$ onto an open subgroup of~$W$.

As a consequence, there is a uniquely determined
topology on~$S'$ making it a topological group,
and which makes $W$ an open subgroup
and induces on it the given locally compact topology.
Throughout the following, $S'$ will be equipped with
this locally compact topology,
which is finer than the
topology induced by~$S$.
Given $p\in F$, let $\pr_p$
be the canonical projection of $\prod_{q\in F}G_q$
onto $G_p$. Pick a compact, open subgroup
$C$ of~$W$.
Considerations very similar
to the preceding ones show
that there is a uniquely
determined group topology on
$H_p:=\pr_p(S')$ such that $\pr_p(C)$,
equipped with the compact topology induced by $G_p$
(which makes it a $p$-adic Lie group)
is a compact, open subgroup of~$H_p$.
Then $H_p$ is a $p$-adic Lie group,
and $\pr_p|_C^{\spr_p(C)}$ being continuous
and open, we see that $\pr_p|_{S'}^{H_p}\!: S'\to H_p$
is a quotient map. The topology on $C$
being induced by the maps $\pr_p|_C^{\spr_p(C)}$,
where $\pr_p(C)$ is open
in $H_p$, we easily see
that $\prod_{p\in F}H_p$ induces the given
locally compact topology on~$S'$.
Note that $D$, being discrete in $S$,
is {\em a fortiori\/} a discrete (and hence closed)
normal subgroup of $S'$ (whose topology is finer).
Hence $X:=\phi(G)\isom S'/D$ is an $\A_\cp$-group,
where we equip $X$ now with
the topology making $\rho|_{S'}^X$
a quotient morphism.
In order that $G$
be an $\A_\cp$-group,
it only remains to show that $\theta:=\phi|^X\!: G\to X$
is an isomorphism of topological groups.
But
$\rho|_W^{\rho(W)}=\rho|_W^{\phi(U)}$
is a quotient morphism with respect to the
new topologies on domain and range.
Since the topology on the domain~$W$
coincides with the old topology,
we deduce that so does the topology on the image $\phi(U)$.
Now $\phi|_U^{\phi(U)}$
being an isomorphism
and $\phi(U)=\rho(W)$ being open in~$X$,
we see that~$\theta$ is an isomorphism.
\end{proof}
\begin{rem}\label{bet1}
Suppose that $G$ is an $\A_\cp$-group,
say $G=S/D$ as in Corollary~\ref{additinfo}.
Applying the construction from
the proof of Proposition~\ref{admissible}\,(b)
to $\phi:=\id_G$, we see that $G=S'/D$
where $S'$ is a closed subgroup
of a product $P:=\prod_{p\in F}H_p$
of $p$-adic Lie groups
for some finite subset $F\sub \cp$,
$D$ is a discrete normal
subgroup of~$S'$, and furthermore
all of the coordinate projections
$\pr_p\!: P\to H_p$ restrict
to {\em quotient morphisms\/}
$\pr_p|_{S'}\!: S'\to H_p$.
\end{rem}
\begin{rem}\label{better}
Let $\cp$ be a finite
set of primes, $S$ a closed subgroup
of a product $\prod_{p\in \cp}G_p$
of $p$-adic Lie groups, and $f\!: G\to S$
be a continuous homomorphism from a locally compact group
to~$S$.
Repeating the proof of
Proposition~\ref{admissible}\,(b)
with $D:=\{1\}$, we see that
$G/{\ker f}$ is isomorphic
to a closed subgroup of a product
$\prod_{p\in \cp}H_p$ of $p$-adic Lie groups.
\end{rem}
\begin{rem}\label{charlsub}
Proposition~\ref{admissible}
allows us to apply Proposition~\ref{appradm}
to the cases $\A:=\A_\cp$ and $\A:=\Sub_\cp$.
We deduce, in particular,
that every $\Mix_\cp$-group (resp., $\LSub_\cp$-group)
is a pro-$\A_\cp$-group (resp., a pro-$\Sub_\cp$-group),
whence it is a
projective limit of a projective system of $\A_\cp$-groups
(resp., $\Sub_\cp$-groups),
such that all bonding maps and all limit maps
are quotient morphisms.
We also deduce the useful fact
that a locally compact group is a
$\Mix_\cp$-group (resp., a $\LSub_\cp$-group)
if and only if it can be approximated by
$\A_\cp$-groups (resp., by $\Sub_\cp$-groups).
\end{rem}
\begin{thm}\label{intersections}
For any sets of primes $\cp$ and $\cq$,
we have
\begin{itemize}
\item[\rm (a)]
$\,\LSub_\cp \cap \LSub_\cq = \LSub_{\cp\cap\cq}\;${\rm ;}
\item[\rm (b)]
$\,\LSub_\cp\cap \Mix_\cq \,= \,\Mix_{\cp\cap \cq}\;${\rm ;}
\item[\rm (c)]
$\,\Mix_\cp\cap \Mix_\cq \;= \;\Mix_{\cp\cap \cq}\;$.
\end{itemize}
\end{thm}
\begin{proof}
We may assume that $\cp,\cq\not=\emptyset$, the excluded
case being trivial.

(a) Let $G\in \LSub_\cp\cap \LSub_\cq$.
By Remark~\ref{charlsub},
in order that $G\in \LSub_{\cp\cap \cq}$,
we only need to show that $G$ can be approximated
by $\Sub_{\cp\cap \cq}$-groups. To verify the latter, let
$U$ be an identity neighbourhood
of~$G$; after shrinking~$U$,
we may assume that $U$ is a compact, open subgroup
of~$G$.
Since $G\in \LSub_\cp$,
there exists a closed normal subgroup
$K\sub U$ of~$G$
such that $G/K\in \Sub_\cp$ (Remark~\ref{charlsub}),
and thus $G/K\in \Sub_F$ for some finite subset $F\sub \cp$.
Now $G\in \LSub_\cq$ entails
that $G/K\in \LSub_\cq$,
whence the identity neighbourhood
$U/K$ contains
a closed normal subgroup~$N$
of $G/K$ such that $(G/K)/N$ has an
open subgroup of the form
$H=\prod_{q\in E}H_q$,
where $H_q$ is a $q$-adic Lie group
for $q$ in a finite subset $E\sub \cq$.
The class $\Sub_F$ being closed under
the formation of Hausdorff quotients
and closed subgroups,
we see that $H$ is a $\Sub_F$-group.
Hence $H$ has an open subgroup
of the form $W=\prod_{p\in F}W_p$
for certain $p$-adic Lie groups
$W_p$;
we may assume that~$W_p$
is a pro-$p$-group.
Given $p\in F$, for each $q\in E \setminus \{p\}$,
the continuous homomorphism
$\pr_q|_{W_p}\!: W_p\to H_q$
(where $\pr_q\!: H\to H_q$ is the coordinate projection)
has kernel~$W_p$,
by Lemma~\ref{propproq}.
Hence $\pr_q(W_p)=\{1\}$
for each $q\in E\setminus \{p\}$.
Then $W_p=\{1\}$ for all $p\in F\setminus E$.
Consequently,
$W=\prod_{p\in E\cap F}W_p$,
showing that $H$ and thus also
$(G/K)/N$ is a $\Sub_{\cp\cap \cq}$-group.
Let $\rho\!: G\to G/K$ be the quotient map.
Then $\rho^{-1}(N)\sub U$ and this
is a closed normal subgroup of~$G$ such that
$G/ \rho^{-1}(N)\isom (G/K)/N\in \Sub_{\cp\cap\cq}$.
Thus $G$ can be approximated by $\Sub_{\cp\cap \cq}$-groups.

(b) Suppose that $G\in \LSub_\cp\cap
\Mix_\cq$.
By Remark~\ref{charlsub},
in order that $G\in \Mix_{\cp\cap \cq}$,
we only need to show that $G$ can be approximated
by $\A_{\cp\cap \cq}$-groups. To verify this, let
$U$ be a compact, open subgroup
of~$G$.
Since $G\in \LSub_\cp$,
there exists a closed normal subgroup
$K\sub U$ of~$G$
such that $G/K\in \Sub_F$
for some finite, non-empty subset $F\sub \cp$.
Now $G\in \Mix_\cq$ entails
that $G/K\in \Mix_\cq$,
whence the identity neighbourhood
$U/K$ contains
a closed normal subgroup~$N$
of $G/K$ such that $(G/K)/N$ is an
$\A_\cq$-group (Remark~\ref{charlsub}).
Hence, in view of Remark~\ref{bet1},
there are $q$-adic Lie groups
$H_q$ for $q$ in a
non-empty finite subset $E\sub \cq$,
a closed subgroup $S\sub \prod_{q\in E}H_q=:H$,
and a quotient morphism $\rho\!: S\to (G/K)/N$,
with discrete kernel~$D$,
such that $\pr_q|_S\!: S\to H_q$
is a quotient morphism
for each~$q$, where $\pr_q\!: H\to H_q$
is the coordinate projection.
Now, $D$ being discrete,
the groups $S$ and $(G/K)/N$ are locally
isomorphic.
The group $(G/K)/N$
is a $\Sub_F$-group as a quotient of
the $\Sub_F$-group $G/K$
(Proposition~\ref{admissible}\,(a)),
whence every identity neighbourhood
of $(G/K)/N$ contains a compact, open subgroup
which is a product of $p$-adic Lie groups,
with $p\in F$. Hence also~$S$ has a compact, open subgroup
of the form $V=\prod_{p\in F}V_p$,
where $V_p$ is a $p$-adic Lie group
for each $p\in F$. Let $q\in E\setminus F$.
By Lemma~\ref{propproq}, $\pr_q|_{V_p}$
has open kernel for each $p\in F$,
whence $\pr_q|_V$ has open kernel.
The group~$V$ being compact,
$\ker\pr_q|_V$ has finite index in $V$,
entailing that $\pr_q(V)$ is finite.
The latter set being open in~$H_q$
(since $\pr_q|_S$ is a quotient morphism),
we deduce that $H_q$ is a discrete
group (and hence a $p$-adic Lie group
for each $p$). Thus
$H=D'\times \prod_{q\in E\cap F}H_q$,
where $D':=\prod_{q\in E\setminus F}H_q$ is discrete.
Therefore $(G/K)/N\isom S/D$ is an $\A_{E\cap F}$-group
and hence an $\A_{\cp\cap \cq}$-group.
The kernel of the natural quotient map
$G\to (G/K)/N$ being contained in~$U$,
we see that $G$ can be approximated by $\A_{\cp\cap \cq}$-groups,
as required.

(c) Since $\Mix_\cp\sub \LSub_\cp$,
assertion\,(c) is a trivial consequence of~(b).
\end{proof}
Note that
$\Mix_{\{p\}}=\LSub_{\{p\}}$,
because a locally compact group
having a $p$-adic Lie group
as an open subgroup
is itself a $p$-adic Lie group.
We now show that $\Mix_\cp$ is a proper
subclass of $\LSub_\cp$ whenever
$\cp$ has at least two elements.
\begin{prop}\label{discern}
If $\cp$ is a non-empty, non-singleton
set of primes, then there exists
a group $G\in \Sub_\cp\sub \LSub_\cp$
such that $G\not\in \Mix_\cp$.
\end{prop}
\begin{proof}
Let $p,q\in \cp$ be two distinct
primes and $n\in \N$
such that both $\SL_n(\Q_p)$
and $\SL_n(\Q_q)$ are simple
groups.
Set $G_1:=\SL_n(\Q_p)\times \SL_n(\Z_q)$,
$G_2:=\SL_n(\Z_p)\times \SL_n(\Q_q)$
and $H:=\SL_n(\Z_p)\times \SL_n(\Z_q)$.
Then the amalgamated product
$G:=G_1*_H G_2$ can be made a topological
group with $H$ as a compact, open
subgroup, and thus $G\in \Sub_{\{p,q\}}$.
Then $G\not\in \Mix_\cp$.
In fact, otherwise $G\in \Mix_\cp\cap
\LSub_{\{p,q\}}=\Mix_{\{p,q\}}$,
whence $\cp=\{p,q\}$ without loss of generality.
To derive a contradiction,
let $U_1\sub \SL_n(\Z_p)$
and $U_2\sub \SL_n(\Z_q)$
be compact, open subgroups
such that $U_1$ is a pro-$p$-group
and $U_2$ a pro-$q$-group.
For simplicity of notation,
we identify $H$, $G_1$ and $G_2$ with the corresponding
subgroups of~$G$,
and we identify $\SL_n(\Q_p)$ with
$\SL_n(\Q_p)\times\{\one\}\sub G_1$
and $\SL_n(\Q_q)$ with
$\{\one\}\times \SL_n(\Q_q)\sub G_2$.
As we suppose that $G\in \Mix_{\{p,q\}}$,
there exists a closed normal subgroup
$N\sub G$ such that $N\sub U_1\times U_2$
and $G/N\in \A_{\{p,q\}}$
(cf.\ Proposition~\ref{admissible}\,(b)).
Since $\SL_n(\Q_p)\cap N$
is a proper normal subgroup
of $\SL_n(\Q_p)$, we must have
$\SL_n(\Q_p)\cap N=\{\one\}$.
Likewise,
$\SL_n(\Q_q)\cap N=\{\one\}$
and thus
$N=(N\cap U_1)\times (N\cap U_2)=\{\one\}$,
exploiting Corollary~\ref{nicesub}.
Hence $G\in \A_{\{p,q\}}$,
and thus there exists a $p$-adic Lie group
$H_1$, a $q$-adic Lie group $H_2$,
a closed subgroup $S\sub H_1\times H_2$
and a closed normal subgroup
$M\sub S$ such that $G\isom S/M$.
Let $\pi\colon S\to G$ be a quotient
morphism with kernel~$M$.
By Corollary~\ref{nicesub},
$(S\cap H_1)\times (S\cap H_2)$
is open in~$S$.
The map~$\pi$ being continuous,
we find compact, open subgroups
$V_1\sub S\cap H_1$
and $V_2\sub S\cap H_2$
such that
$V_1$ is a pro-$p$-group,
$V_2$ is a pro-$q$-group,
and $\pi(V_1\times V_2)\sub U_1\times U_2$.
Then $\pi(V_1)\sub U_1$ and $\pi(V_2)\sub U_2$,
by Lemma~\ref{propproq}.
Since $\pi(V_1)$ is open in $\SL_n(\Q_p)$,
we find $x_1\in V_1\sub H_1$
such that $\pi(x_1)$ is not a diagonal
matrix and has, say,
a non-zero $(i,j)$-entry
(where $i\not=j$).
Let $d$
be the diagonal matrix
whose $i$-th and $j$-th diagonal
entries are $p^{-k}$ and $p^k$, respectively,
while all other diagonal entries are~$1$.
Choosing $k\in \N$ large enough, we
obtain $z_1:=d \pi(x_1)d^{-1}\not\in \SL_n(\Z_p)$.\
There is $g\in S$ such that $\pi(g)=d$.
Then $h_1:=gx_1g^{-1}\in S\cap H_1$
(the latter subgroup being normal in~$S$),
and $\pi(h_1)=z_1\in \SL_n(\Q_p)$. 
Likewise,
we find an element $h_2\in S\cap H_2$
such that $z_2:=\pi(h_2)\in \SL_n(\Q_q)$
but $z_2\not\in \SL_n(\Z_q)$.
Then $h_1h_2=(h_1,h_2)=h_2h_1$.
Thus $h_1$ and $h_2$ commute and hence
also their images~$z_1$ and~$z_2$
under~$\pi$ commute.
But $z_1z_2\not=z_2z_1$.
To see this, for $j\in \{1,2\}$
choose a set $R_j\sub G_j$
of representatives
for the cosets in $G_j/H$,
such that $z_j\in R_j$.
Then the word $z_1z_2$
in $G_1*G_2$
is the normal form of $z_1z_2$
(as in \cite[Thm.\,11.66]{Rot}),
and the word $z_2z_1$
is the normal form of $z_2z_1$.
Hence $z_1z_2\not=z_2z_1$ in~$G$,
and we have reached
the desired contradiction.
\end{proof}
\section{Tidy subgroups and the scale function for
{\boldmath $\Sub_\cp$\/}-groups and {\boldmath $\LSub_\cp$\/}-groups}
We now describe tidy subgroups
and calculate the scale function
for $\Sub_\cp$-groups
and then, passing to projective limits,
for $\LSub_\cp$-groups.
For this purpose, we
require a
slight generalization of the notion of the module
of an automorphism.
\begin{numba}
If $G$ is a locally compact group,
with Haar measure $\lambda$, and $\alpha\!: H\to G$
an injective, continuous homomorphism from an open subgroup
$H\sub G$ onto
an open subgroup $S:=\alpha(H)$ of~$G$,
then, due to uniqueness of Haar measure on~$S$
up to a multiplicative constant,
there is a positive
real number $\Delta_G(\alpha)$
(also written $\Delta(\alpha)$ when $G$ is understood),
the {\em module of~$\alpha$},
such that
$\lambda|_S = \Delta_G(\alpha) \, \alpha(\lambda|_H)$.
Thus $\Delta_G(\alpha)=\frac{\lambda(\alpha(U))}{\lambda(U)}$,
for every non-empty open subset $U\sub H$
of finite measure.
If $G$ is a $p$-adic Lie group here,
identifying $L(H)$ and $L(S)$ with $L(G)$, we have
\begin{equation}\label{genBou}
\Delta_G(\alpha)\; =\; \Delta_{L(G)}(L(\alpha))\; =\;
|\dt\, L(\alpha)|_p\,,
\end{equation}
using the natural
absolute value $|.|_p$ on $\Q_p$.
In fact, the proof of
\cite[Ch.\,3, \S3.16, Prop.\,55]{BLi}
(treating only \'{e}tale endomorphisms)
directly generalizes to the present situation.
\end{numba}
It is also useful to know
that every $\Sub_\Primes$-group
has an open subgroup
satisfying the ascending chain condition
on closed subgroups,
because this property
ensures that a compact, open subgroup
satisfying condition (T1)
of tidiness (as described in the Introduction)
automatically satisfies
condition (T2) as well
\cite[Thm.\,3.32 and Rem.\,3.33\,(2)]{BaW}.
\begin{prop}\label{chain}
Every $G\in \Sub_\Primes$
has a compact, open subgroup~$U$
such that $U$ satisfies the
ascending chain condition on closed subgroups.
\end{prop}
\begin{proof}
Since $G\in \Sub_\Primes$,
there exists a finite set $F\sub \Primes$
and an open subgroup $U$ of~$G$
such that $U=\prod_{p\in F} U_p$
for certain $p$-adic Lie groups~$U_p$.
After shrinking~$U_p$,
we may assume that $U_p$ also
is a pro-$p$-group.
We may furthermore assume that each~$U_p$
satisfies the ascending chain condition
on closed subgroups, because
every $p$-adic
Lie group has an
open subgroup
with this property~\cite[proof of Prop.\,3.5]{Wan}.
Now let $S_1\sub S_2\sub \cdots$
be an ascending sequence
of closed subgroups of~$U$.
Then $(S_n\cap U_p)_{n\in \N}$
becomes stationary for each $p\in F$,
and hence so does $(S_n)_{n\in \N}$,
as $S_n=\prod_{p\in F} (S_n\cap U_p)$
by Proposition~\ref{block}.
\end{proof}
Note that $\Z_p^\N$ is a $\Mix_\Primes$-group
without open subgroups satisfying an ascending chain condition
on closed subgroups.\\[3mm]
The following fact is essential
for the calculation of the scale function and tidy subgroups.
\begin{numba}\label{defncontr}
If $E$ is a finite-dimensional
$\Q_p$-vector space
and $\alpha$
a linear automorphism of~$E$,
then
\begin{eqnarray}
\!\!  E\, \;  &= & E_{\rm p}\oplus E_0\oplus E_{\rm m}\,, \quad\mbox{where}
\label{contrdeco}\\
E_{\rm m} &:=& \{x\in E\!: \mbox{$\alpha^n(x)\to 0$ as $n\to\infty$}\}\,,
\nonumber\\
E_{\rm p}\,
 &:=& \{x\in E\!: \mbox{$\alpha^{-n}(x)\to 0$ as $n\to\infty$}\}
\quad\mbox{and}\nonumber \\
\!\! E_0\,  &:=& \{x\in E\!: \mbox{$\alpha^\Z(x)$ is relatively compact}\}\,;
\nonumber
\end{eqnarray}
see \cite[La.\,3.4]{Wan} or \cite[La.\,3.3]{Sca}
(cf.\ \cite[pp.\,80--83]{Mar}
for more refined information).
We call (\ref{contrdeco})
the {\em contraction decomposition of~$E$
with respect to~$\alpha$.} 
It is known
(see, e.g., \cite[La.\,3.3]{Sca} and its proof)
that
there is an ultrametric norm $\|.\|$ on~$E$
which is {\em adapted\/}
to the decomposition~(\ref{contrdeco})
in the sense that
$\|\alpha(x)\|=\|x\|$ for all $x\in E_0$
and, for suitable $\theta>1$,
\[
\|\alpha(x)\|\; \geq \; \theta \, \|x\|\quad
\mbox{for all $\,x\in E_{\rm p}$}\quad
\mbox{and}\quad
\|\alpha(x)\|\; \leq \; \theta^{-1}\|x\|\quad
\mbox{for all $\, x\in E_{\rm m}$}\,.
\]
\end{numba}
If $\g$ is a finite-dimensional
$p$-adic Lie algebra, there exists
a compact, open submodule
$V\sub \g$ such that the Campbell-Hausdorff
series converges on $V\times V$
to a function $*\!: V\times V\to V$
making~$V$ a $p$-adic Lie group
(see \cite[Ch.\,II, \S8.3, Prop.\,3]{BLi}).
We call $V$ a {\em CH-group\/}.
\begin{thm}\label{scaletidy}
Let $\cp\not=\emptyset$ be a finite set of primes,
$G$ be a $\Sub_\cp$-group,
and $\alpha\!: G\to G$ be an automorphism
$($e.g., $\alpha:=I_x\!: y\mto xyx^{-1}$
for some $x\in G)$.
Let $H$ be an open subgroup of~$G$
of the form $\prod_{p\in \cp}H_p$,
where $H_p$ is a $p$-adic Lie group
for $p\in \cp$. Then we have:
\begin{itemize}
\item[\rm (a)]
For each $p\in \cp$,
there exists an  open subgroup~$U_p$ of $H_p$
such that $\alpha(U_p)\sub H_p$.
\item[\rm (b)]
Identifying $L(U_p)$ with $L(H_p)$,
we consider $\beta_p:=L (\alpha|_{U_p}^{H_p})$ as a Lie algebra
automorphism of $L(H_p)$.
We let $L(H_p)=L(H_p)_{\rm p}\oplus L(H_p)_0\oplus L(H_p)_{\rm m}$
be the contraction decomposition of $L(H_p)$
with respect to $\beta_p$, and
abbreviate
$L(H_p)_+:=L(H_p)_{\rm p}\oplus L(H_p)_0$
and $L(H_p)_-:=L(H_p)_0\oplus L(H_p)_{\rm m}$.
Then
\begin{equation}\label{explscale}
r_G(\alpha)\;
=\;\prod_{p\in \cp} \;\Delta \bigl(\beta_p|_{L(H_p)_+}^{L(H_p)_+}\bigr)
\;=\;
\prod_{p\in \cp} \;\; \prod_{\stackrel{\scriptstyle
i\in \{1,\ldots,\dim L(H_p)\}}{\scriptstyle|\lambda_{p,i}|_p\geq 1}}
\!|\lambda_{p,i}|_p\;,
\end{equation}
where $\lambda_{p,1},\ldots ,\lambda_{p,\dim L(H_p)}$
are the eigenvalues of $\beta_p$ in an algebraic
closure $\wb{\Q_p}$ of $\Q_p$
$($repeated according to their algebraic
multiplicities$)$, and $|.|_p$ is the unique
extension of the usual absolute value
on $\Q_p$ to an absolute value on $\wb{\Q}_p$.
In particular,
\begin{equation}\label{imagesca}
\im\, s_G \; \sub \; \prod_{p\in \cp}p^{\N_0}\,,\quad
\mbox{i.e.,}\quad  \Primes(G)\; \sub \; \cp\, .\vspace{-2.8mm} 
\end{equation}
\item[\rm (c)]
For each~$p$, let $\phi_p\!: V_p\to W_p$
be a topological isomorphism from a CH-group
$V_p\sub L(H_p)$ onto an open subgroup $W_p$ of~$H_p$,
such that $L(\phi_p)=\id_{L(H_p)}$.
Let $\|.\|_p\!: L(H_p)\to [0,\infty[$
be an ultrametric norm on $L(H_p)$ adapted to the
contraction decomposition.
Then there exists $\ve_0>0$
such that
$B_{p,\ve_0}:=\{y\in L(H_p)\!: \|y\|_p<\ve_0\}\sub V_p$
for each $p\in \cp$, and such that
\begin{equation}\label{tidyexpl}
B_\ve \; :=\; \prod_{p\in \cp}\phi_p(B_{p,\ve})
\end{equation}
is a compact, open subgroup of~$G$ which is tidy for $\alpha$,
for each $\ve\in \;]0,\ve_0]$.
\end{itemize}
\end{thm}
\begin{proof}
Let $H=\prod_{p\in \cp}H_p$ be as described in the theorem.
After shrinking~$H_p$,
we may assume that $H_p$
is a pro-$p$-group and that
there exists an
isomorphism $\phi_p\!: V_p\to H_p$
with $L(\phi_p)=\id_{L(H_p)}$,
for some compact, open submodule $V_p\sub L(H_p)$,
equipped with the Campbell-Hausdorff
multiplication.
There exists a compact, open subgroup
$W_p\sub H_p$ such that $\alpha(W_p)\sub H$
and $\alpha^{-1}(W_p)\sub H$.
Then
$\pr_q\circ \alpha|_{W_p}^H=1$
and $\pr_q\circ \alpha^{-1}|_{W_p}^H=1$
for $q\not=p$,
by Lemma~\ref{propproq},
where $\pr_q\!: H\to H_q$ is the coordinate projection.
Hence
\begin{equation}\label{goodW}
\alpha(W_p)\sub H_p\quad\mbox{and}\quad
\alpha^{-1}(W_p)\sub H_p\,,\quad\mbox{for each $p\in \cp$.}
\end{equation}
Shrinking $W_p$ further, we may assume that
$\alpha|_{W_p}^{H_p}$ is linear in exponential coordinates,
viz.\
\begin{equation}\label{loclin}
\phi_p(L(\alpha|_{W_p}^{H_p}).y)=\alpha(\phi_p(y))\qquad
\mbox{for all $\,y\in \phi_p^{-1}(W_p)$,}
\end{equation}
and likewise for $\alpha^{-1}$.
In the following, we identify $H_p$ with $V_p\sub L(H_p)$
by means of the isomorphism $\phi_p^{-1}$,
for convenience. Then (\ref{loclin}) and
its analogue for $\alpha^{-1}$ take the form
\begin{equation}\label{simplerf}
\alpha(y)=L(\alpha|_{W_p}^{H_p}).y\quad\mbox{and}\quad
\alpha^{-1}(y)=L(\alpha|_{W_p}^{H_p})^{-1}.y\quad \mbox{for all $y\in W_p$.}
\end{equation}
For each $p\in \cp$, we choose
an ultrametric norm $\|.\|_p$ on $L(H_p)$
adapted to the contraction decomposition
$L(H_p)=L(H_p)_{\rm p}\oplus L(H_p)_0\oplus L(H_p)_{\rm m}$
of $L(H_p)$
with respect to the Lie algebra
automorphism $\beta_p:=L(\alpha|_{W_p}^{H_p})$.
Choose $\delta>0$ such that
$B_{p,\delta}:=\{y\in L(H_p)\!: \|y\|_p< \delta\}\sub H_p$
for each $p\in \cp$. Next, choose
$\ve_0\in \; ]0,\delta]$ such that
\[
B_{p,\ve_0}\sub W_p\quad \mbox{and}\quad  \alpha(B_{p,\ve_0})\sub
B_{p,\delta},\quad\mbox{for all $p\in \cp$,}
\]
and such that $B_{p,\ve}$ is a subgroup of $(H_p,*)$
(and hence of~$G$),
for all $\ve\in \;]0,\ve_0]$;
the latter is possible by
\cite[Ch.\,III, \S4.2, La.\,3\,(iii)]{BLi}.
For $p\in \cp$, consider the map
\[
f_p\!:
(W_p\cap L(H_p)_+)\times (W_p\cap L(H_p)_{\rm m})\to H_p\sub L(H_p),\quad
(u,v)\mto uv
\]
(product in~$G$).
Then $f_p$ is analytic
(and hence strictly differentiable by \cite[4.2.3 \& 3.2.4]{BVa})
and its
differential at $0$ is the identity map $\id_{L(H_p)}$.
By the Inverse Function Theorem
in the form
\cite[Prop.\,7.1\,(b)$'$]{IMP},
after shrinking $\ve_0$ we can achieve
that
\begin{equation}\label{nicetrick}
f_p\bigl(
(B_{p,\ve}\cap L(H_p)_+)\times (B_{p,\ve}\cap L(H_p)_{\rm m})
\bigr)
\; =\; B_{p,\ve}\, ,\quad\;
\mbox{for each $\, \ve\in \; ]0,\ve_0]$.}
\end{equation}
Fix $\ve\in \; ]0,\ve_0]$;
we claim that the compact, open subgroup $B_\ve:=\prod_{p\in \cp}B_{p,\ve}$
of~$G$ is tidy for~$\alpha$.
To this end, note that
\[
\alpha^{-1}(B_{p,\ve}\cap L(H_p)_+)
\; =\;
\beta_p^{-1}.(B_{p,\ve}\cap L(H_p)_+)
\; \sub \;
B_{p,\ve}\cap L(H_p)_+
\]
by choice of the norm $\|.\|_p$, and thus
$\alpha^{-n}(B_{p,\ve}\cap L(H_p)_+)\sub B_{p,\ve}$
for each $n\in \N_0$,
entailing that
$B_{p,\ve}\cap L(H_p)_+\sub (B_\ve)_+$
and thus
\[
\prod_{p\in \cp} (B_{p,\ve}\cap L(H_p)_+)\; \sub \; (B_\ve)_+,
\]
where $(B_\ve)_\pm:=\bigcap_{n\in \N_0}\alpha^{\pm n}(B_\ve)$.
If $y=(y_p)_{p\in \cp}\in B_\ve$
and $y_q\not \in L(H_q)_+$ for some $q\in \cp$,
then $\|(\beta_q)^{-n}.y_q\|_q
\geq \|(\beta_q)^{-n}.(y_q)_{\rm m}\|_q\geq \theta^n\,\|(y_q)_{\rm m}\|_q
\to \infty$ as $n\to\infty$,
where $(y_q)_{\rm m}$ is the component of $y_q$ in
$L(H_q)_{\rm m}$ and $\theta>1$
is as in {\bf\ref{defncontr}},
the definition of an adapted norm.
This entails that there is
$m\in \N$ such that $(\beta_p)^{-n}(y_p)\in B_{p,\ve}$
for all $n\in \{0,\ldots,m-1\}$
and all $p\in \cp$,
and an element $q\in \cp$
such that
\begin{equation}\label{outside}
(\beta_q)^{-m}.y_q \; \not\in \; B_{q,\ve}\,.
\end{equation}
Then $\beta_p^{-n}.y_p=\alpha^{-n}(y_p)$ for
all $n\in \{1,\ldots, m\}$ and all $p\in \cp$,
by~(\ref{simplerf}).
Since
$\alpha^{-m+1}(y_p)\in B_{p,\ve}$,
we have $\alpha^{-m}(y_p)\in H_p$ for each $p\in \cp$.
Hence (\ref{outside})
entails that $\alpha^{-m}(y)=(\beta_p^{-m}.y_p)_{p\in \cp}\not\in B_\ve$,
whence $y\not\in (B_\ve)_+$.
Summing up, we have shown that
\[
(B_\ve)_+ \; =\; \prod_{p\in \cp} (B_{p,\ve}\cap L(H_p)_+)\,,
\]
and an analogous argument gives $(B_\ve)_-=
\prod_{p\in \cp} (B_{p,\ve}\cap L(H_p)_-)$.
Using
(\ref{nicetrick}),
we see that
\begin{eqnarray*}
B_\ve & \supseteq & (B_\ve)_+(B_\ve)_- \; =\; \prod_{p\in \cp}
(B_{p,\ve}\cap L(H_p)_+)(B_{p,\ve}\cap L(H_p)_-)\\
&\supseteq&
\prod_{p\in \cp}
(B_{p,\ve}\cap L(H_p)_+)(B_{p,\ve}\cap L(H_p)_{\rm m})
\; =\; \prod_{p\in \cp}B_{p,\ve}\; =\; B_\ve
\end{eqnarray*}
and hence $B_\ve=(B_\ve)_+(B_\ve)_-$,
{\em i.e.},
$B_\ve$ satisfies condition (T1)
of tidiness.
As a consequence of Proposition~\ref{chain},
$B_\ve$ also satisfies (T2)
and thus $B_\ve$ is tidy for~$\alpha$.

To calculate $r_G(\alpha)$, we choose $\ve\in \; ]0,\ve_0]$
and obtain
\begin{eqnarray*}
r_G(\alpha) &=& [\alpha((B_\ve)_+):(B_\ve)_+]
\; =\; \prod_{p\in \cp}\,
[\alpha(B_{p,\ve}\cap L(H_p)_+)\,:\, B_{p,\ve}\cap L(H_p)_+]\\
&=&\prod_{p\in \cp}\Delta(\alpha|_{B_{p,\ve}\cap L(H_p)_+}^{H_p\cap L(H_p)_+})
\;
\stackrel{{\rm (5)}}{=}
\;
\prod_{p\in \cp} \Delta (\beta_p|_{L(H_p)_+})
\; =\;  \prod_{p\in \cp} |\dt (\beta_p|_{L(H_p)_+})|_p\\
&=&
\prod_{p\in \cp} \;\; \prod_{\stackrel{\scriptstyle
i\in \{1,\ldots,\dim L(H_p)\}}{\scriptstyle|\lambda_{p,i}|_p\geq 1}}
\!|\lambda_{p,i}|_p\;,
\end{eqnarray*}
with $\lambda_{p,i}$ as described in the theorem.
Here, we used \cite[La.\,3.4]{Sca}
to pass to the third line.
Since $|\Q_p^\times|_p=p^{\N_0}$,
the second term in the second line
shows that $r_G(\alpha)\in \prod_{p\in \cp}p^{\N_0}$.
Hence $\,\im\,s_G\sub \prod_{p\in \cp}p^{\N_0}$
in particular.\vspace{-1mm}
\end{proof}
\begin{rem}
In the special case where~$G$ is a $p$-adic
Lie group, Theorem~\ref{scaletidy}
provides a self-contained,
explicit calculation of the scale
function, and a basis of subgroups
tidy for $x\in G$.
The earlier calculation of $s_G$ in~\cite{Sca}
relied on a result from~\cite{Wan},
and tidy subgroups could not be described
explicitly in that paper.
\end{rem}
\begin{la}\label{furth}
Let $G$ be a totally disconnected,
locally compact group,
$K\sub G$ be a compact, normal subgroup,
$q\!: G\to G/K$ be the quotient map,
and $x\in G$. Then $q^{-1}(U)$ is tidy for $x$,
for every compact, open subgroup $U\sub G/K$
which is tidy for~$xK$, and $s_G(x)=s_{G/K}(xK)$.
\end{la}
\begin{proof}
Let $y:=xK$ and $V:=q^{-1}(U)$.
Since $I_x^n(V)$
is $K$-saturated for each $n\in\Z$,
and $q^{-1}(I_y^n(U))=I^n_x(V)$,
we easily see that $V_\pm=q^{-1}(U_\pm)$,
$V=q^{-1}(U_+U_-)=V_+V_-$, and $V_{++}=q^{-1}(U_{++})$,
which is closed. Thus $V$ is
tidy for~$x$, and $s_G(x)=[I_x(V_+):V_+]=[q^{-1}(I_y(U_+)):q^{-1}(U_+)]
=[I_y(U_+):U_+]=s_{G/K}(y)$.
\end{proof}
Combining Theorem~\ref{scaletidy}
and Lemma~\ref{furth}, we obtain:
\begin{cor}\label{scalelim}
Let $\cp\not=\emptyset$ be a set of primes,
and
$G\in \LSub_\cp$.
Then we have:
\begin{itemize}
\item[\rm (a)]
Let $K\sub G$ be a compact, normal subgroup
such that $G/K$ is a $\Sub_\cp$-group;
then $G/K$ is a $\Sub_F$-group
for some finite subset $F\sub \cp$.
We have
\begin{equation}
s_G(x)=s_{G/K}(xK)\qquad \mbox{for each $x\in G$,}
\end{equation}
where $s_{G/K}(xK)$ can be calculated explicitly as
described in Theorem~{\rm \ref{scaletidy}}.
In particular, $\Primes(G)\sub F$,
whence $\Primes(G)$ is a finite subset of~$\cp$.
\item[\rm (b)]
Let $x\in G$.
For every $K$ as in {\rm (a)}
and compact, open subgroup $U\sub G/K$
tidy for $xK$, the subgroup $q_K^{-1}(U)\sub G$
is tidy for~$x$, where $q_K\!: G\to G/K$
is the quotient map. The set
of subgroups $q_K^{-1}(U)$ tidy for~$x$,
for all possible~$K$ and~$U$
as before,
is a basis for the filter
of identity neighbourhoods of~$G$.
Hence $G$ has small tidy subgroups.\Punkt
\end{itemize}
\end{cor}
\begin{rem}\label{adele}
If $G$ is an Ad\`{e}le group then
$\Primes(G/{G_0})$ typically is an infinite set.
In this case,
$G/{G_0}$ is not a $\Mix_\Primes$-group
(nor a $\LSub_\Primes$-group),
by Corollary~\ref{scalelim}\,(a).
For example, we have
$\Primes(G/G_0)=\Primes$ for $G:=\dl_F
\bigl(\SL_n(\R)\times\prod_{p\in F}
\SL_n(\Q_p)\times \prod_{\Primes\setminus F}\SL_n(\Z_p)\bigr)$,
where $F$ ranges through the set of finite subsets
of~$\Primes$, and $n\geq 2$ (cf.\ \cite[Thm.\,5.1]{Sca}).
\end{rem}
\section{The minimal set of primes needed to build up
a compactly generated $\Mix_\Primes$-group, or $\LSub_\Primes$-group}
As a tool, we introduce an analogue of
the adjoint action of Lie groups
for $\Sub_\Primes$-groups.
\begin{numba}\label{defnAdp}
Suppose that $G$ is a $\Sub_\Primes$-group
and $H\sub G$ an open subgroup of the form
$H=\prod_{p\in F}H_p$,
where $F$ is a finite set of primes
and $H_p$ is a $p$-adic Lie group, for each
$p\in F$. Given $x\in G$, consider the inner automorphism
$I_x\!: G\to G$, $I_x(y):=xyx^{-1}$.
Given $x$,
for each $p\in F$, there is a compact, open subgroup
$U_p\sub H_p$ such that $I_x(U_p)\sub H_p$.
Identifying $L(U_p)$ with $L(H_p)$ by means of
the isomorphism of Lie algebras $L(i_p)$,
where $i_p\!: U_p\to H_p$ is the inclusion map,
we may consider $\Ad_p(x):=L(I_x|_{U_p}^{H_p})$
as an automorphism of the Lie algebra $L(H_p)$.
It is easy to see that
\[
\Ad_p\!: G\to \Aut(L(H_p)),\quad x\mto \Ad_p(x)
\]
is a homomorphism.
Since $\Ad_p(x)=\Ad(\pr_p(x))$
for $x\in H$, where $\pr_p\!: H\to H_p$ is
the coordinate projection and $\Ad\!: H_p\to \Aut(L(H_p))$
is continuous, we see that the homomorphism
$\Ad_p$ is continuous on the open subgroup $H$ and hence
continuous.
\end{numba}
\begin{thm}\label{exactprimes}
Let $\cp$ be a set of primes.
\begin{itemize}
\item[\rm (a)]
If $G\in \LSub_\Primes$
is compactly generated,
then $G\in \LSub_\cp$
if and only if $\,\Primes(G)\sub \cp$.
\item[\rm (b)]
If $G\in \Mix_\Primes$
is compactly generated,
then $G\in \Mix_\cp$
if and only if $\,\Primes(G)\sub \cp$.
\end{itemize}
\end{thm}
\begin{proof}
By Corollary~\ref{scalelim},
$G\in \LSub_\cp$
entails $\Primes(G)\sub \cp$
for $\cp$ non-empty,
and apparently $\Primes(G)\sub \cp$ also holds if
$\cp=\emptyset$,
as every pro-discrete group is uniscalar.

(a) Let $G\in \LSub_\Primes$
be compactly
generated; we want to show that $G\in \LSub_{\Primes(G)}$.
In view of Remark~\ref{charlsub},
we only need to show that $G$ can be approximated by
$\Sub_{\Primes(G)}$-groups.
Now, given a compact, open subgroup $U\sub G$,
there exists a compact normal subgroup
$N\sub U$ of~$G$ such that
$G/N$ is a $\Sub_\cq$-group
for some
finite subset $\cq\sub \Primes$.
Then $G/N$ is compactly
generated, and
$\Primes(G)=\Primes(G/N)$ (Lemma~\ref{furth}),
whence $\Primes(G)\sub \cq$, by what has already been
shown.
Let $\rho\!:G\to G/N$
be the quotient map.
If we can show that $G/N\in \LSub_{\Primes(G)}$,
then we can find a compact normal subgroup
$Z\sub \rho (U)$ of $G/N$ such that $(G/N)/Z\isom G/\rho^{-1}(Z)$
is a $\Sub_{\Primes(G)}$-group, where
$\rho^{-1}(Z)\sub U$, whence indeed $G$ can be approximated
by $\Sub_{\Primes(G)}$-groups.

Replacing $G$ with $G/N$, we may
therefore assume that $G\in \Sub_\cq$
for some finite set of primes~$\cq$.
Let $U\sub G$
be as before.
Then $G$ has an open subgroup $H\sub U$ of the form
$H=\prod_{p\in \cq}H_p$,
where $H_p$ is a $p$-adic Lie group.
Given $x\in G$,
define $I_x\!: G\to G$, $I_x(y):=xyx^{-1}$.
After shrinking $H_p$, we may identify
$H_p$ with a compact, open $\Z_p$-submodule
of $L(H_p)$, equipped with the CH-multiplication
(as in the proof of Theorem~\ref{scaletidy}\,(a))
and may assume that $I_x(y)=\Ad(x).y$
for all $x,y\in H_p$.
We let $\emptyset\not=K$ be a compact, symmetric
generating set for~$G$.
Then $K\sub FH$
for some finite subset $F\sub K$.
For each $p\in \cq$ and $x\in F$,
there exists a compact, open subgroup
$V_p(x)\sub H_p$ such that $I_x(V_p(x))\sub H_p$
and $I_x(y)=L(I_x|_{V_p(x)}^{H_p}).y$
for all $y\in V_p(x)$;
we set $V_p:=\bigcap_{x\in F}V_p(x)$.
Let $W_p$ be a compact, open, normal
subgroup of~$H_p$ such that $W_p\sub V_p$.
Since $K\sub FH$, where $I_x|_{H_q}\ident 1$
for $x\in H_p$ with $p\not=q\in \cq$,
we conclude that
$I_x(W_p)\sub H_p$
for all $x\in K$, and
$I_x(y)=L(I_x|_{W_p}^{H_p}).y=
\Ad_p(x).y$ for all $x\in K$ and $y\in W_p$.
Now let $\Cr:=\cq\setminus\Primes(G)$.
Given $p\in \Cr$ and $x\in G$
we deduce from
(\ref{explscale}) and the fact that $p$ neither divides
$s_G(x)$ nor $s_G(x^{-1})$ that all eigenvalues
of $L(I_x|_{W_p}^{H_p})$ in $\wb{\Q_p}$
have modulus~1.
Repeating the arguments
used to prove ``1)$\impl$4)'' of Prop.\,3.1
in \cite{GaW}, we find that $\Ad_p(x)=L(I_x|_{W_p}^{H_p})$
is a compact element of $\Aut(L(H_p))$,
for each $x\in G$ and each $p\in \Cr$.
The homomorphism $\Ad_p\!: G\to \Aut(L(H_p))$
being continuous (see {\bf \ref{defnAdp}}),
we deduce that the subgroup
$R_p:=\Ad_p(G)\sub \Aut(L(H_p))$
is compactly generated.
Being compactly generated and periodic, $R_p$
is relatively compact in $\Aut(L(H_p))$
(see \cite{Par}).
As a consequence of \cite[Part~II, Appendix~1, Thm.\,1]{Ser},
there exists a compact, open $\Z_p$-submodule
$M_p\sub W_p$ of $L(H_p)$
which is invariant under $R_p$.
Then the subgroup $C_p:=\langle M_p\rangle$
of $W_p$ generated by~$M_p$
is open in $W_p$ and compact,
and it is a normal subgroup of~$G$ as
it is normalized by each $x\in K$,
where
$K$ generates~$G$
(here we use that
$I_x(C_p)=\langle I_x(M_p)\rangle=\langle
\Ad_p(x).M_p\rangle=\langle M_p\rangle=C_p$).
As a consequence,
$C := \prod_{p\in \Cr}C_p \sub U$
is a compact, normal subgroup
of~$G$ such that
$G/C$ contains
\[
H/C\; \isom \; \prod_{p\in \Primes(G)}H_p\; \times\;
\prod_{p\in \Cr}H_p/C_p
\]
as an open subgroup,
where $\prod_{p\in \Cr}H_p/C_q$
is discrete. Thus $G/C\in \Sub_{\Primes(G)}$,
showing that $G$ can be approximated by $\Sub_{\Primes(G)}$-groups,
which completes the proof of~(a).

(b) Now suppose that $G\in \Mix_\Primes$
is compactly
generated.
Then $G\in \LSub_\Primes$ {\em a fortiori\/}
and hence $G\in \LSub_{\Primes(G)}$, by Part\,(a).
Thus $G\in \Mix_\Primes \cap \LSub_{\Primes(G)}=
\Mix_{\Primes\cap \Primes(G)}=\Mix_{\Primes(G)}$,
using Theorem~\ref{intersections}\,(b).
\end{proof}
Generalizing \cite[Cor.\,5]{Par} and \cite[Thm.\,5.2]{GaW},
we obtain:
\begin{cor}\label{unispro}
Every compactly generated, uniscalar $\Mix_\Primes$-group
$($or $\LSub_\Primes$-group$)$ $G$
is pro-discrete.
\end{cor}
\begin{proof}
Since $\Primes(G)=\emptyset$,
Theorem~\ref{exactprimes} shows that
$G\in \LSub_\emptyset$.
\end{proof}
Theorem~\ref{exactprimes} and Corollary~\ref{unispro}
become false for groups that are not compactly generated,
as there is a
uniscalar $p$-adic Lie group without a compact,\hspace*{-.3mm}
open,\hspace*{-.3mm} normal
subgroup~\cite[\S\,6]{GaW}.
\section{Variants based on locally pro-{\boldmath $p$\/}
groups}\label{seclocpro}
Variants of some of our
results can be obtained when
$p$-adic Lie groups are replaced
with locally pro-$p$ groups.
We also provide counterexamples
for
results which do not
carry over.
\begin{numba}
In {\bf\ref{deflocpro}},
we introduced the class $\Loc_p$
of locally pro-$p$ groups.
Given
$\emptyset\not= \cp\sub \Primes$,
we set $\Loc_\cp:=\bigcup_{p\in \cp}\Loc_p$.
Then
$\cV(\Loc_\cp)={\rm SC}(\A_\cp^\vee)$,
where
$\A_\cp^\vee := {\rm \wb{Q}\,\wb{S}P}(\Loc_\cp)$.
We define
\[
\Mix_\cp^\vee\;
:=\; \{G\in \cV(\Loc_\cp)\!: \,\mbox{$G$ is locally compact}\,\}\,,
\]
$\Mix_\emptyset^\vee:=\Mix_\emptyset$,
and $\A_\emptyset^\vee:=\A_\emptyset$.
Given $\cp\sub \Primes$,
we let $\Sub_\cp^\vee$
be the class of all topological groups
possessing an open subgroup
isomorphic to $\prod_{p\in F}G_p$,
where $F\sub \cp$ is finite
and $G_p$ a pro-$p$-group
for each $p\in F$.
We define
\[
\LSub_\cp^\vee \; :=\;
\{G\in \cV(\Sub_\cp^\vee)\!: \,\mbox{$G$ is locally compact}\,\}\,.
\]
In particular,
$\LSub_\emptyset^\vee=\LSub_\emptyset$.
\end{numba}
\begin{numba}\label{moregen}
Then
$\A_\cp\sub \A_\cp^\vee$
and $\Sub_\cp\sub \Sub_\cp^\vee$ and thus
$\Mix_\cp\sub \Mix_\cp^\vee$ and
$\LSub_\cp\sub \LSub_\cp^\vee$.
\end{numba}
\begin{thm}\label{portmanteau}
All of Corollary~{\rm\ref{injthensub}}--Theorem~{\rm \ref{intersections}}
remain valid if $p$-adic Lie groups
are replaced with locally pro-$p$ groups,
$\A_\cp$ by $\A_\cp^\vee$,
$\Mix_\cp$ by $\Mix_\cp^\vee$,
$\Sub_\cp$ by $\Sub_\cp^\vee$,
and $\LSub_\cp$ by $\LSub_\cp^\vee$.
With
analogous replacements,
also Corollary~{\rm\ref{scalelim}\,(a)} carries
over.\,\footnote{Except for the reference to
Theorem~\ref{scaletidy} for the explicit
calculation of $s_{G/K}(xK)$.}
\end{thm}
\begin{proof}
The proofs of Corollary~\ref{injthensub}--Theorem~\ref{intersections}
can be repeated verbatim in the new situation,
making the replacements
described in the theorem.
Also the adaptation of
Corollary~\ref{scalelim}\,(a)
is immediate
in view of the Proposition~\ref{scalep} below.\vspace{-.6mm}
\end{proof}
\begin{rem}\label{convrem}
In particular, $\A_\cp^\vee$ and $\Sub_\cp^\vee$
are suitable for approximation.
Hence a locally compact group
belongs to $\Mix_\cp^\vee$
and $\LSub_\cp^\vee$ if and only if it can be approximated
by $\A_\cp^\vee$-groups (resp.,
by $\Sub_\cp^\vee$-groups).
\end{rem}
\begin{prop}\label{scalep}
Let $G$ be a locally pro-$p$ group.
Then
\[
r_G(\alpha)\; \in \; p^{\N_0}
\]
for each $\alpha\in \Aut(G)$
and thus $\Primes(G)\sub \{p\}$, i.e.,
the scale function $s_G$ takes its values in $p^{\N_0}$.
In particular, the conclusions
apply if
$G$ is an analytic Lie group $($or $C^1$-Lie group$)$
over a local field $\K$,
and $p:=\car(\ck)$ the characteristic
of the residue field $\ck$ of~$\K$.
\end{prop}
\begin{proof}
Let $U$ be an open, pro-$p$ subgroup of~$G$.
Since~$U$ cannot contain an
infinite $q$-Sylow subgroup
for any $q\not=p$,
\cite[end of p.\,173]{Wi2}
entails that $r_G(\alpha)\in p^{\N_0}$
(see also Section~\ref{seclocal}).
Since every $C^1$-Lie group
over~$\K$ is locally pro-$p$
(see {\bf \ref{deflocpro}}),
the final assertion follows.
\end{proof}
\begin{rem}
Of course,
\cite[Thm.\,3.5]{Sca} (and Theorem~\ref{scaletidy}
above) provide much more
information in
the special
case where $G$ is a Lie group over a local
field of characteristic~$0$.
\end{rem}
\begin{rem}
If $G$ is a Lie group over local
field of positive characteristic,
then not every topological group
automorphism of $G$ needs to be
analytic.
It is interesting that
Proposition~\ref{scalep}
provides some information also
on these non-analytic automorphisms.
The scale of analytic automorphisms
has been studied in~\cite{POS}.
While
\begin{equation}\label{char0case}
r_G(\alpha)\;=\; r_{L(G)}(L(\alpha))
\end{equation}
always holds in characteristic~$0$,
surprisingly this equation becomes
false in general if $\car(\K)>0$.
Closer inspection reveals
that~(\ref{char0case})
holds if and only if
$G$ has
small subgroups tidy for~$\alpha$
(see \cite{POS}).
This
natural 
property (which is not always satisfied)
was first explored in~\cite{BaW}
in the context of contraction groups
and has been exploited further in~\cite{TID}.
\end{rem}
A straightforward adaptation
of the proof of Theorem~\ref{intersections}\,(a) also shows:
\begin{prop}\label{overlap}
$\,\LSub_\cp\cap \LSub_\cq^\vee=\LSub_{\cp\cap\cq}$
holds,
for all subsets $\,\cp,\cq \sub \Primes$.\Punkt
\end{prop}
Here are some differences.
First of all,
$\Sub_\Primes^\vee$-groups
need not have an
open subgroup satisfying an
ascending chain condition
on closed subgroups,
as the pro-$p$ group $\Z_p^\N$ shows.
Next, Theorem~\ref{exactprimes}
and Corollary~\ref{unispro}
do not carry over to $\Mix_\cp^\vee$-groups
and $\LSub_\cp^\vee$-groups.
\begin{example}
If we start with a non-trivial,
finite $p$-group~$K$,
the construction described
in \cite[p.\,269, lines 2--5]{BMP}
(based on an ansatz from~\cite{KaW})
outputs a totally disconnected, locally compact group~$G$
which is compactly generated and
uniscalar, but does not possess a compact open
normal subgroup. Since, by construction,
$G$ contains an open subgroup topologically
isomorphic to $K^\N$, we see that $G$ is
locally pro-$p$.
\end{example}
Finally, we observe that
$\Mix_{\{p\}}^\vee$-groups need not be
locally pro-$p$.
Indeed:
Given any prime $q\not=p$
and non-trivial finite $q$-group~$K$,
the group $G:=K^\N$
is pro-discrete and
thus $G\in \Mix_{\{p\}}^\vee$.
However, $G$ is not locally pro-$p$.
\section{The minimal set of primes in the general case}\label{seclocal}
In this section,
we show that also
for groups $G\in \Mix_\Primes$
that are not compactly
generated, there always is a smallest
set of primes $\cp$ such that
$G\in \Mix_\cp$.
As a tool to find~$\cp$,
we associate certain sets
of primes to totally
disconnected, locally compact
groups~$G$, which
only depend on the local isomorphism type
of~$G$.
\begin{defn}\label{content}
Given a totally disconnected,
locally compact group~$G$,
we let $\bL(G)\sub \Primes$
be the set of all primes $p$
such that, for
every compact, open subgroup
$U\sub G$,
the element $p$ divides
the index $[U:V]$
for some compact, open subgroup
$V\sub U$.
The set $\bL(G)$ is called the
{\em local prime content\/}
of~$G$.
The {\em reduced prime content\/} of~$G$
is defined as
\[
\bL_r(G)\; :=\; {\bigcap}_K \; \bL(G/K)\,,\vspace{-2mm}
\]
where $K$ runs through the set
of all
compact, normal subgroups
of~$G$. Then $\bL_r(G)\sub \bL(G)$.
\end{defn}
Standard arguments from
Sylow theory show that
$p\in \bL(G)$
if and only if some (and hence any)
compact, open subgroup~$U$
of~$G$ has an infinite
$p$-Sylow subgroup.
In \cite[end of p.\,173]{Wi2},
it has been noted that this is the case
if $p\in \Primes(G)$.
More generally, we observe:
\begin{prop}\label{asymptotic}
Let $G$ be a totally disconnected,
locally compact group. Then we have:
\begin{itemize}
\item[\rm (a)]
For every $\alpha\in \Aut(G)$,
the set of prime divisors of $\, r_G(\alpha)$
is a subset of $\, \bL(G)$.
\item[\rm (b)]
$\,\Primes(G)\sub \bL(G)$ holds,
and indeed $\, \Primes(G)\sub \bL_r(G)$.
\end{itemize}
\end{prop}
\begin{proof}
(a) If $p\in \Primes\setminus \bL(G)$,
then there exists a compact, open subgroup
$U$ of~$G$ such that $p$ and $[U:V]$
are coprime, for every compact, open subgroup
$V\sub U$.
Now $U_N:=\bigcap_{n=0}^N \alpha^n(U)$
satisfies condition (T1)
of tidiness for~$\alpha$, for some $N\in \N_0$
(see \cite[La.\,1]{Wi1}).
By
\cite[La.\,2.2]{Wi2},
$r_G(\alpha)$ divides $[U_N: U_N\cap \alpha^{-1}(U_N)]$, which in turn
divides $[U:U_N\cap \alpha^{-1}(U_N)]$,
because
$[U:U_N\cap \alpha^{-1}(U_N)]=
[U:U_N]\cdot [U_N:U_N\cap \alpha^{-1}(U_N)]$.
Since $p$ and $[U,U_N\cap\alpha^{-1}(U_N)]$
are coprime, we deduce that $p$ does not divide
$r_G(\alpha)$.

(b) Is immediate from (a) and Lemma~\ref{furth}.
\end{proof}
We shall use the following simple observations:
\begin{la}\label{opers}
\begin{itemize}
\item[\rm (a)]
If $G$ is a totally disconnected,
locally compact group
and $U\sub G$ an open subgroup,
then $\bL(G)=\bL(U)$.
\item[\rm (b)]
If $G:=\prod_{j=1}^n G_j$
is a finite product of
totally disconnected, locally
compact groups $G_j$, then
$\bL(G)=\bigcup_{j=1}^n \bL(G_j)$.
\end{itemize}
\end{la}
\begin{proof}
The proof is obvious
from the definition
of the local prime contents,
the definition of the product
topology
and the fact that $[U:W]=[U:V]\cdot
[V:W]$ for any compact, open subgroups
$W\sub V\sub U$.
\end{proof}
Proposition~\ref{scalep} can be rephrased as follows:
\begin{prop}\label{openprop}
Let $G$ be a locally pro-$p$ group
$($e.g., a $C^1$-Lie group over
a local field whose residue field
has characteristic~$p)$.
Then $\bL(G)\sub \{p\}$,
and $\bL(G)=\{p\}$
if and only if $G$ is non-discrete.
In particular,
$\Primes(G)\sub \{p\}$.
\end{prop}
\begin{proof}
Let $U\sub G$ be a compact, open subgroup
which is pro-$p$.
Then $\bL(G)=\bL(U)\sub \{p\}$,
where $\bL(U)=\emptyset$ if and only if~$U$
is discrete. The rest follows from
Proposition~\ref{asymptotic}\,(b).
\end{proof}
Our next aim is to analyze
$\LSub_\Primes^\vee$-groups by means of their
``intermediate'' prime content:
\begin{defn}\label{defninterm}
If $G$ is a $\LSub_\Primes^\vee$-group,
we define its {\em intermediate prime content\/}
via
\[
\bL_i(G)\; :=\;
\bigcup_{K'\in J(G)}\;\bigcap_{K'\supseteq K\in J(G)} \bL(G/K)\,,\vspace{-2mm}
\]
where $J(G)$ is the set of all
compact, normal subgroups $K$ of~$G$
such that $G/K\in \Sub_\Primes^\vee$.
\end{defn}
Note that $J(G)$ is a filter basis in~$G$
converging to~$1$, by Proposition~\ref{appradm}\,(b)
and Remark~\ref{convrem}.\\[3mm]
The following lemma explains the terminology ``intermediate'':
\begin{la}\label{specialK}
If $G$ is a $\LSub_\Primes^\vee$-group, then
\begin{equation}\label{spec}
\bL_r(G)= \bigcap_{K\in J(G)} \bL(G/K)\,,\vspace{-2mm}
\end{equation}
where $J(G)$ is as in Definition~{\rm \ref{defninterm}}.
Furthermore, $\,\bL_r(G)\sub \bL_i(G)\sub \bL(G)$.
\end{la}
\begin{proof}
The inclusion ``$\sub$'' in (\ref{spec})
is obvious. To see the converse inclusion,
let $p\in \Primes\setminus \bL_r(G)$.
Then
there exists a compact, normal subgroup
$K$ of~$G$ such that $p\not\in \bL(G/K)$.
Hence, there exists a compact, open subgroup
$U\sub G/K$ such that $p$ and the index
$[U:V]$ are coprime, for every compact, open subgroup
$V$ of~$U$. Let
$\pi\!: G\to G/K$ be the quotient map.
Since
$G/K$ is a $\LSub_\Primes^\vee$-group,
there exists a compact, normal
subgroup $N$ of $G/K$
such that
$H:=(G/K)/N$
is a $\Sub_\Primes^\vee$-group
(Remark~\ref{convrem}).
Let $q\!: G/K\to H$
be the quotient morphism.
Then $U':=q(U)$
is a compact, open subgroup
of $H$ such that
$[U':V]=[q^{-1}(U'):q^{-1}(V)]=[U:q^{-1}(V)]$
is not divisible by~$p$,
for any compact, open subgroup
$V$ of $H$, and thus $p\not\in \bL(H)$.
Since $G/M\isom H\in \Sub_\Primes^\vee$,
where $M:=\pi^{-1}(N)$
is a compact normal subgroup
of~$G$, we see that $p$
is not contained in the right hand side of (\ref{spec}).
Thus (\ref{spec}) is established.

It is obvious
that $\bL_r(G)\sub \bigcap_{K'\supseteq K\in J(G)}
\bL(G/K)$ for all $K'\in J(G)$, and thus $\bL_r(G)\sub \bL_i(G)$.
To see that $\bL_i(G)\sub \bL(G)$,
let $p\in \bL_i(G)$.
Then there exists $K'\in J(G)$
such that $p\in \bigcap_{K'\supseteq K\in J(G)}\bL(G/K)$.
For every compact, open subgroup
$U\sub G$, there exists $K''\in J(G)$
such that $K''\sub U$ (see Remark~\ref{convrem}).
Hence, $J(G)$ being a filter basis,
we find $K\in J(G)$ such that $K\sub K'\cap K''\sub U$.
Let $\rho\!: G\to G/K$ be the quotient map.
Since $p\in \bL(G/K)$,
there exists a compact, open subgroup
$V\sub \rho(U)$ such that $[\rho(U):V]=[U:\rho^{-1}(V)]$
is divisible by~$p$.
The subgroup $\rho^{-1}(V)$
of~$U$ being compact and open,
we deduce that $p\in \bL(G)$.\vspace{-3mm}
\end{proof}
\begin{thm}\label{noncomp}
Let $\cp$ be a set of primes.
Then the following holds:
\begin{itemize}
\item[\rm (a)]
If $\,G\in \Sub_\cp^\vee$,
then $\,\bL(G)\sub \cp$.
\item[\rm (b)]
If $\,G\in \LSub_\Primes^\vee$,
then $G\in \LSub_\cp^\vee$ if and only if
$\,\bL_i(G)\sub \cp$.
\item[\rm (c)]
If $\,G\in \Mix_\Primes^\vee$,
then $G\in \Mix_\cp^\vee$
if and only if $\,\bL_i(G)\sub \cp$.
\item[\rm (d)]
If $\,G\in \LSub_\Primes$,
then $G\in \LSub_\cp$ if and only if
$\,\bL_i(G)\sub \cp$.
\item[\rm (e)]
If $\,G\in \Mix_\Primes$,
then $G\in \Mix_\cp$
if and only if $\,\bL_i(G)\sub \cp$.
\end{itemize}
\end{thm}
\begin{proof}
(a)
There is a finite subset $F\sub \cp$ such that
$G$ has an open subgroup
of the form $U=\prod_{p\in \cp}U_p$,
where $U_p$ is a non-discrete pro-$p$-group
for $p\in F$.
Hence $\bL(G) = F\sub \cp$,
by Lemma~\ref{opers}
(a), (b) and Proposition~\ref{openprop}.

(b) If $G\in \LSub_\cp^\vee$,
then $\bL(G/K)\sub\cp$
for each $K\in J(G)$, by (a).
Hence $\bL_i(G)\sub \cp$.
Let us show now that $G\in \LSub_{\bL_i(G)}^\vee$,
for each $G\in\LSub_\Primes^\vee$. It suffices to show
that $G$ can be approximated
by $\Sub_{\bL_i(G)}^\vee$-groups.
Thus, let $U\sub G$ be a compact, open subgroup
of~$G$. There exists a compact, normal subgroup
$K'\sub U$ of~$G$ such that
$G/{K'}\in \Sub_\Primes^\vee$.
Hence, there exists a finite set of primes
$F$ such that $G/{K'}$ has an open subgroup
of the form $V=\prod_{p\in F}V_p$,
where $V_p$ is a non-discrete pro-$p$-group
for each $p\in F$.
We claim that $F\sub \bL_i(G)$ (whence
indeed $G$ can be approximated by $\Sub_{\bL_i(G)}^\vee$-groups).
To see this, suppose to the contrary
that there exists some $\bar{p}\in F\setminus \bL_i(G)$.
Then $\bar{p}\not\in \bigcap_{K'\supseteq K\in J(G)}\bL(G/K)$,
whence there exists $K\in J(G)$ such that
$K\sub K'$ and $\bar{p}\not\in \bL(G/K)$.
Let $\rho\!: G/K\to G/{K'}$ be the natural map,
$\rho(gK):=gK'$.
Then $G/K$ has an open
subgroup
$W\sub \rho^{-1}(V)$ of the form
$W=\prod_{p\in E}W_p$ for some finite set
of primes~$E$, where $W_p$ is a non-discrete
pro-$p$-group for each $p\in E$.
Since $\bL(G/K)=E$
by the proof of~(a),
we see that $\bar{p}\not\in E$,
whence $\pr_{\bar{p}}\circ
\rho|_{W_p}^V=1$
for each $p\in E$
by Lemma~\ref{propproq},
where $\pr_{\bar{p}}\!: V\to V_{\bar{p}}$ is the coordinate
projection.
Hence
$\pr_{\wb{p}}(\rho(W))=\{1\}$.
The latter being an open subset
of~$V_{\bar{p}}$, we deduce that $V_{\bar{p}}$
is discrete. We have reached a contradiction.

(c)--(e): If $G\in \Mix_\Primes^\vee$,
then $G\in \Mix_\cp^\vee$ if and only if
$G\in \LSub_\cp^\vee$ (by the analogue of
Theorem~\ref{intersections}\,(b)
subsumed by Theorem~\ref{portmanteau}).
Similarly, Proposition~\ref{overlap}
entails that $G\in \LSub_\Primes$ belongs
to $\LSub_\cp$
if and only if $G\in \LSub_\cp^\vee$,
and hence $G\in \Mix_\Primes$ belongs to
$\Mix_\cp$ if and only if
$G\in \LSub_\cp^\vee$, by Theorem~\ref{intersections}\,(b).
Thus (d)--(e)
follow from~(b).\vspace{-.6mm}
\end{proof}
\begin{rem}
The argument used to prove
Part\,(b) of the preceding theorem
shows that $\bL(G/K)\supseteq \bL(G/{K'})$
for all $K,K'\in J(G)$ such that $K\sub K'$.
Hence $\bL_i(G)$ is in fact given by the
simpler formula
$\bL_i(G)=\bigcup_{K\in J(G)}\bL(G/K)$,
for each
$\LSub_\Primes^\vee$-group~$G$.
\end{rem}
\begin{rem}\label{thesmallest}
Given $G\in \LSub_\Primes^\vee$,
Theorem~\ref{noncomp}\,(a)
shows
that $\cp:=\bL_i(G)$
is the smallest set
of primes such that $G\in \LSub_\cp^\vee$.
Parts\,(b)--(e)
have analogous interpretations.
Strengthening Theorem~\ref{intersections}\,(a), (c)
and its analogue
in Theorem~\ref{portmanteau}, we deduce that\vspace{-1mm}
\[
\bigcap_{i\in I} \; \LSub_{\cp_i}^\vee\;\, =\;\,  \LSub_\cp^\vee\qquad
\mbox{with $\, \cp:=\bigcap_{i\in I} \cp_i\,$,}\vspace{-3mm}
\]
for any family $(\cp_i)_{i\in I}$ of sets $\cp_i\sub \Primes$.
Analogous formulas hold for $\LSub_\cp$, $\Mix_\cp$
and $\Mix_\cp^\vee$.
\end{rem}
\begin{rem}
If $G\in \LSub_\Primes$,
instead of $J(G)$ we can
use the set $\wb{J}(G)$
of all compact, normal subgroups $K\sub G$
such that $G/K\in \Sub_\Primes$
to define a set $\wb{\bL}_i(G)$ analogous to $\bL_i(G)$.
Repeating the preceding proofs
with $\LSub_\Primes$ instead of $\LSub_\Primes^\vee$,
we see that $\cp:=\wb{\bL}_i(G)$
is the smallest set of primes
with $G\in \LSub_\cp$ and thus $\wb{\bL}_i(G)=\bL_i(G)$.
\end{rem}
$\LSub_\emptyset^\vee$ being the class of locally compact, pro-discrete
groups, Theorem~\ref{noncomp}\,(b) implies:
\begin{cor}
A group $G\in \LSub_\Primes^\vee$ is
pro-discrete if and only if $\bL_i(G)=\emptyset$.\Punkt
\end{cor}
The following corollary is immediate
from Theorem \ref{exactprimes} and Theorem~\ref{noncomp}\,(d):
\begin{cor}
If $G$ is a compactly generated
$\LSub_\Primes$-group, then
$\Primes(G)=\bL_i(G)$.\Punkt
\end{cor}
\appendix
\section{The set of normal subgroups
with Lie quotients need not be a filter basis}
Let $\K=\R$ or $\K=\Q_p$
for some~$p$.
Given a topological group~$G$, let
$N_\K(G)$
be the set
of all closed normal subgroups $N\sub G$
such that $G/N$ is a $\K$-Lie group.
We describe
a complete abelian topological group~$G$
such that $N_\K(G)$ is not a filter basis.
Examples for such behaviour
had not been known before.
For locally compact~$G$,
the pathology cannot occur,
the class of $\K$-Lie groups
being suitable for approximation
(see \cite[1.7]{App}; cf.\ \cite{HMS}).\\[3mm]
{\bf Construction of {\boldmath $G$}.}
The topology induced by $\K$ on $\Q$
can be refined to a topology~$\tau$ which makes
$\Q$ a non-discrete, complete topological group~\cite{Mrn}.
We write $H:=(\Q,\tau)$
and define $G:=\K\times H$.
Then $H$ is not a $\K$-Lie group, as it is
countable but non-discrete.
Hence $G$ is not a $\K$-Lie group either, since otherwise
$H\isom G/(\K\times\{0\})$ would be a $\K$-Lie group.
The first
coordinate projection
$G\to \K$, $(x,y)\mto x$
is a quotient homomorphism,
with kernel $N:=\{0\}\times H$.
On the other hand, the inclusion map
$\iota\!: H\to \K$ being continuous,
\[
q\!: G\to \K \, , \quad (x,y)\, \mto \, x+\iota(y)\, =\, x+y
\]
is a continuous homomorphism.
Apparently~$q$ is surjective.
Furthermore, $q$ is open
since, for any $0$-neighbourhoods
$U\sub \K$ and $V\sub H$,
we have $U \sub q(U\times V)$.
We set $M:=\ker q=\{(x,-x)\!: x\in \Q\}$.
Then $M\cap N=\{0\}$, $G$ is not a $\K$-Lie group,
and both $G/M$ and $G/N$ are topologically isomorphic
to~$\K$.\Punkt\\[2.5mm]
{\footnotesize {\em Acknowledgements.}
The research was supported by DFG grant
447 AUS-113/22/0-1
and ARC grant LX 0349209.
The author thanks G.\,A. Willis
for his suggestion to extend the studies
from the $p$-adic case to locally pro-$p$ groups.
His remarks also helped to
find the counterexamples proving
Proposition~\ref{discern}.}\vspace{-3mm}
\noindent
{\footnotesize
{\bf Helge Gl\"{o}ckner}, TU~Darmstadt, FB~Mathematik~AG~5,
Schlossgartenstr.\,7, 64289 Darmstadt, Germany.\\
E-Mail: gloeckner@mathematik.tu-darmstadt.de}
\end{document}